# New electric vehicle charging rate design: An MPEC assessment

Icaro Silvestre Freitas Gomes,[a][1] Adam F. Abdin,[a] Jakob Puchinge,[a,b] Yannick Perez[a]

[a]*Université Paris-Saclay, CentraleSupélec, Laboratoire Génie Industriel, 3 rue Joliot-Curie 91190 Gif-sur-Yvette, France*

[b]*Institut de Recherche Technologique SystemX, Palaiseau 91120, France*

**Abstract**

A high penetration of electric vehicles (EVs) will deeply impact the management of electric power systems. To avoid costly grid reinforcements and the risk of load curtailment due to EV charging, indirect load control via adapted economic signals is a solution proposed by many utilities. Charging costs can be reduced with a domestic tariff applied only to EV charging using a dedicated load measurement method while enhancing the flexibility offered by EVs. We develop a game-theoretical model expressed and treated as a mathematical programme with equilibrium constraints (MPEC) to capture the interaction between a national regulatory authority (NRA) designing these tariffs and heterogeneous agents. First, we analyse the conditions in which EV-only tariffs can be applied for domestic slow charging sessions by comparing different energy profiles. Second, we study the impact of EV charging on different tariff structures to identify the most efficient way of recovering network costs. Submetering with a pure volumetric tariff can bring yearly gains varying from $64 to $110 compared to a flat rate. This depends on the share of investment in grid reinforcement that remains to be made. Finally, we derive policy implications from the results and earmark more sophisticated tariff designs for further investigation.

**Highlights:**

- We build a game-theoretical framework to assess the cost-efficiency of EV-only tariff adoption.[2]
- Full innovators able to invest in DERs with EVs have higher total savings with classic time-of-use rates.
- Capacity charges may increase fairness issues among heterogeneous users.
- Consumers adopting EV-only rates with submetering will perceive yearly gains varying between $64 and $110.

[1] Corresponding author, contactable at: *Université Paris-Saclay, CentraleSupélec, Laboratoire Génie Industriel, 3 rue Joliot-Curie 91190 Gif-sur-Yvette, France*. Email: icarosilvestre4@gmail.com, Tel: +33 7 68 18 59 16.

[2] Abbreviations: BESS, battery energy storage system; EV, electric vehicle; PV, photovoltaics; DER, distributed energy resources; ICE, internal combustion engine.

## 1. Introduction

The shift towards a low-carbon economy requires a great reduction in $CO_2$ emissions coming from the transport sector, which accounts for 24% of direct emissions (IEA, 2020). To reduce $CO_2$, electric vehicles (EVs) are required to substitute ICE vehicles to achieve the decarbonisation goals established during international environmental summits. Between 150 and 230 million vehicles are expected to be on the world's roads by 2030, potentially causing future issues for power systems. This study aims to consider the challenges that grid management might face as a result of the high penetration of EVs (Salah et al., 2015). To avoid costly grid reinforcements and the risk of load curtailment due to EV charging, indirect price control via adapted economic signals is a solution adopted by many utilities (Knezoviˊc et al., 2017). These economic signals, given to EV users via a network tariff and energy price profiles, can have different structures. Users who seek to minimise their electricity bills are then led to a different final utilisation pattern. Most household electricity meters do not separate the rates used for household electricity needs and for charging privately owned EVs. This is known as the 'wholehouse' rate. Today, users have the option of a domestic time-of-use (TOU) rate applied exclusively to EV charging. This tariff is called the EV-only rate. Users of this measurement method can reduce charging costs and enhance the flexibility offered by EVs via adapted price signals.

Several pilots have been conducted in the U.S. (California, Minnesota, Texas) to test the technical feasibility and customer acceptance of these rates (Smart Electric Power Alliance, 2019). Californian electric investor-owned utilities (IOUs) already offer these types of rates in their portfolio for dwellings. For instance, Pacific Gas & Electricity (PG&E) and San Diego Gas & Electric (SD&E) allow residential customers to be billed at a tiered rate for home appliances while for EV charging a specific time-of-use plan is adopted (PGE, 2021, SDGE, 2021). However, as installing a second meter is mandatory for this, EV-only rates have not been widely adopted due to the high associated cost of extra equipment. The need for separate wiring and metering renders such market segmentation expensive and cumbersome (Borenstein et al., 2021). One alternative to avoid upfront second meter costs or fees for residential customers is submetering. In this case, the metering infrastructure inside the electric vehicle supply equipment (EVSE) can be used to measure the electricity coming from the grid specifically used for EV charging. Technological progress in smart meters, communication networks and data management will allow the submetering configuration to be adopted by many utilities for billing purposes. In the US case, an official decision on submetering by the regulator in California is expected in 2021 after the conclusion of submetering pilots (CPUC, 2021).

Advanced metering infrastructure (AMI) for electricity is being rolled out in different places in the

E.U. and U.S. For instance, as of July 2018 all but two member states in the E.U. had conducted at least one cost-benefit analysis of a large-scale rollout of smart meters to at least 80% by 2020. However, only a few of these member states reached the target, while the majority postponed the achievement of this milestone to 2030 (European Commission, 2020). In the U.S., by the end of 2020 75% of U.S. households already had an AMI installation in their homes (Cooper and Shuster, 2021).

To the best of our knowledge, no study has assessed the effects on users and network tariff design of adopting two different rates for the same household simultaneously. We investigate this configuration supposing that the submetering scheme has a dedicated protocol allowing communication between the owner's EVSE and the utility for billing purposes. Our integrated approach considering diverse energy profiles and the impact on tariff design allows a fair investigation of the cost-effectiveness of solutions. We find that a total gain between $ 64 and $ 110 per year is achievable depending on the state of the grid. Users willing to adopt this solution, when possible, can collect the spread between flat and time-of-use profiles, which is not offset by network tariff increases in the case of volumetric tariffs.

The structure of the paper is as follows. First, an overview of the problem is given to explain the motivation for the research with a literature review. The data used are then presented along with the setup proposed. In section 5 the results are presented showing energy profile assessments and network tariff impacts. In section 6 we discuss the results and derive policy implications. The last section concludes.

## 2. Literature review

In this section we analyse two main strands of literature. The first concerns the interaction between EVs, distributed energy resources and tariff design, which has received much attention recently. This is mainly due to the great penetration of DERs in the grid, which may change the way utilities charge their customers. The second looks at EV demand-side flexibility assessments that analyse how smart charging and V2G can simultaneously bring remuneration via energy services and support the grid. This work locates itself at the intersection between these two strands by tackling dedicated EV tariffs which can defer network investment by adopting smart charging.

A vast body of literature investigates the impacts of different electricity rates, including energy prices and network tariffs, on specific end-users possessing DERs. The impacts on their decisions are assessed either with exogenously defined tariffs (Ansarin et al., 2020; Backe et al., 2020; Freitas Gomes et al., 2021; Avau et al., 2021) or using equilibrium models in which grid tariffs are determined endogenously as a result of a bi-level approach (Hoarau and Perez, 2019; Schittekatte and Meeus, 2020; Askeland et al., 2021). In all these cases, only one tariff structure is analysed at a time for users with only one metering scheme. Hoarau and Perez (2018) analyse the interaction between tariff design and DERs with EVs to measure the most cost-efficient configuration and the fairest for heterogeneous agents. For these authors the more a tariff structure gives incentives for DERs, the less beneficial it is for EVs. Using a similar framework, Askeland et al. (2021) highlight that an EV agent can spread EV charging evenly throughout the day to minimise the agents' individual peak load regardless of the overall load situation. This can be problematic if prospective costs drive grid costs since the coincident peak can go up and be followed by a tariff increase. If recovered via cost-reflective tariffs, this type of cost can benefit both prosumers who can invest in DERs and consumers who are not able to, according to Schittekatte and Meeus (2020).

A subset of this literature specifically investigates tariffs and best practices for EV charging (King and Datta, 2018; Hildermeier et al., 2019; Kufeoglu et al., 2019). Hildermeier et al. (2019) argue that customer education is key, in particular to attract new user groups who are not already convinced of the specific advantages of managed EV charging. For EV-only tariffs, this effort would probably duplicate since more information regarding bill savings optimisation and metering infrastructure would be needed. However, if adequate information about tariffs and behaviour is provided, users usually choose more complex tariffs (Mayol and Staropoli, 2021). In the same line, Kufeoglu et al. (2019) support the idea that energy utilities must offer consumers more options for TOU tariffs, not only to allow for greater demand-side management but also to encourage uptake of V2H technology. Finally, King and Datta (2018) point out that submetering is indeed a far less expensive option than installing a separate meter for the EV. Moreover, they explain that standards are already available to ensure appropriate billing using the submeter approach. Its accuracy can be set at the same level as already existing standards for electricity meters.

The second strand of literature studying EV demand-side flexibility raises attention to the challenges and opportunities expected during EV uptake. Knezovic et al. (2017) provide a roadmap with key recommendations for supporting active EV involvement in grids to provide flexibility services such as investment deferral, load and voltage services. One barrier they highlight is a lack of standardised smart-meter functionalities and interoperability among all participants. In the same line, Salah et al. (2015) argue that price incentives can help to exploit available load flexibility embedded in EV charging, while if ill-designed they can lead to a significant increase in peak loads in times of low prices. Last, it is worth pointing out the value streams of numerous V2G services (Thompson and Perez, 2019). In our framework, network investment deferral supported by EV flexibility in a context of renewable energy as explored by Hemmati and Mehrjerdi (2020) can contribute to an accessible cost-reflective tariff.

## 3. Method

In this section the game-theoretic model is presented. First, the solution method to determine the main results is summarised. Then, the equations of the optimisation problems for both levels are described in detail.

*3.1. Model overview*

This section describes the model used to pursue the main goals of this research. We develop a bi-level game-theoretical optimisation model to capture the interaction between a national regulatory authority (NRA) and dwellings. The upper level represents the NRA responsible for setting network tariffs to maximise social welfare subject to grid cost-recovery, total DER investment costs and the cost of final electricity use in dwellings. In the lower level, dwellings seek to minimise their electricity costs according to the final tariff structure applied. Decisions by dwellings include ones on new investments in DER and optimal scheduling of electric vehicle charging, which in turn depends on the tariff applied. Therefore, there is a clear interdependence between lower-level charging, investment decisions and tariff design. This interdependence requires an equilibrium solution to be found, which can only be properly captured with a bi-level optimisation approach as previously described. At first, we only deal with the lower level to find the most suitable energy profile for each agent owning an EV. In this case, the lower level is formulated as mixed integer linear programming (MILP) taking network tariffs as exogenous variables. Then, to assess the impact on network tariffs, a bi-level model is created and treated as a mathematical programme with equilibrium constraints (MPEC) in which the equilibrium game-theoretic solution for both the NRA and dwellings will be the one in which no unilateral deviation in their decisions is profitable, defined as a Nash equilibrium. The complete solution framework is summarised in Fig. 1. The full model formulation and the techniques used to transform and solve the MPEC are detailed step by step in Appendix A.

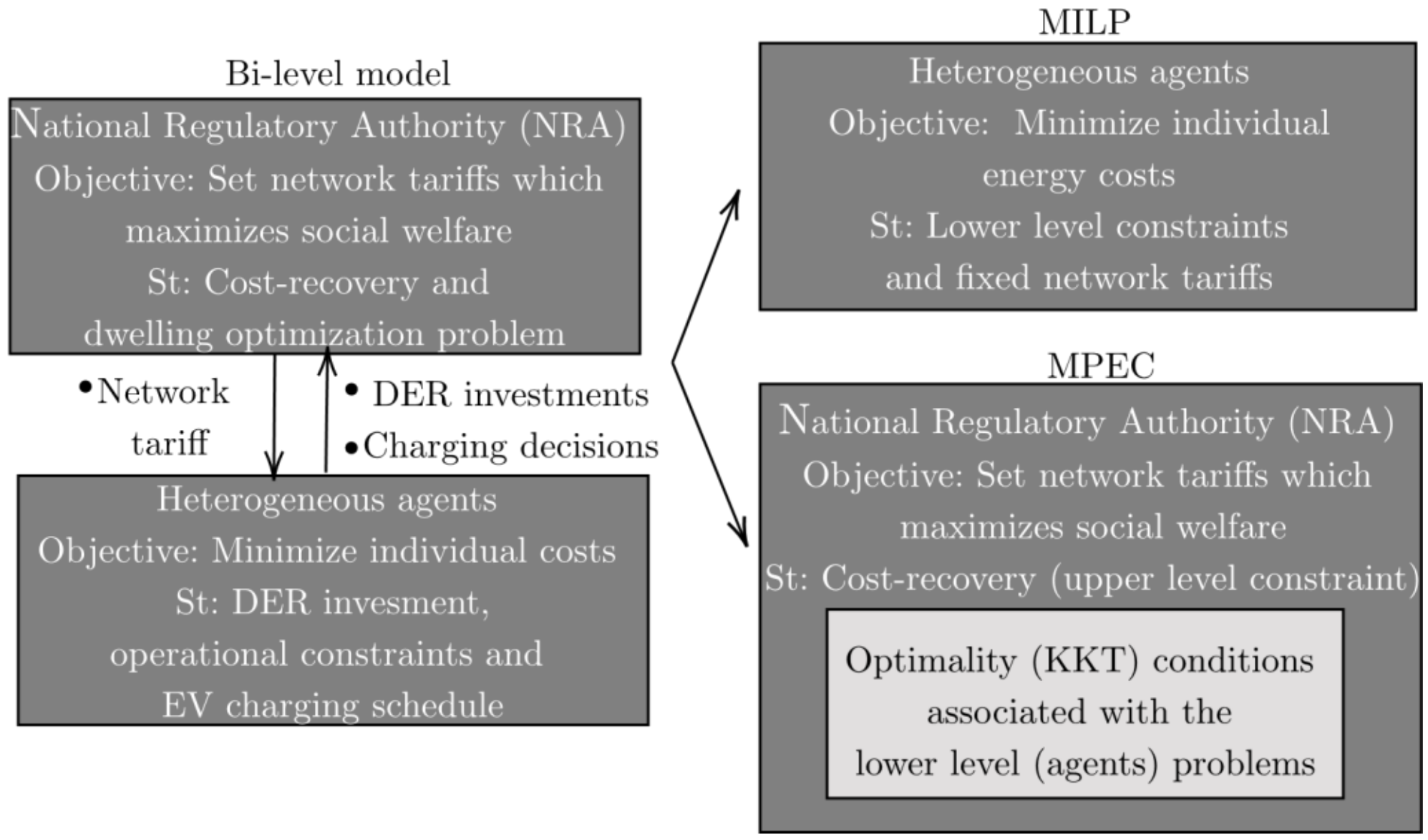

Figure 1: Methodological framework.

*3.2. Upper-level problem formulation: the regulator*

*3.2.1. Regulator objective function*

The main objective of the regulator is to decide on the network tariff to maximise social welfare while ensuring cost-recovery and resolving individual optimisation problems of agents. In this context, social welfare maximisation can be assumed to be equivalent to cost minimisation for all interacting agents since, in addition to network costs, investment and energy costs for dwellings are also included in the regulator's minimisation problem. As a consequence, Eq. (1), which represents the NRA cost function, accounts for the sum of DER investments costs ($Cost^{DER}$), energy charges ($Cost^{P}$) and network charges ($Cost^{N}$) for agents:

$$\min Cost^{NRA} = Cost^{DER} + Cost^{P} + Cost^{N} \tag{1}$$

The first term in the regulator's objective function is investment decisions by dwellings. It is represented by Eq. (2), in which a customer *c* seeks to optimally decide to install certain solar PV ($ic_c^{PV}$) and battery ($ic_c^{S}$) capacities if it is profitable to do so, with annualised investment costs $I^{PV}$ and $I^{S}$ for solar PV and the battery respectively. This is formulated as:

$$Cost^{DER} = \sum_{c \in C} \left( I_c^S * ic_c^S + I_c^{PV} * ic_c^{PV} \right) \tag{2}$$

The total net energy costs to satisfy electricity demand for all agents is calculated in Eq. (3). We assume one retailer supplying all the customers. However, more than one energy profile can be considered according to the rate chosen. For instance, charging vehicles with an EV-only tariff will require a time-of-use profile whereas house electricity consumption can be considered flat. Energy charges account for the total energy imports ($imp_{c,h}^{P}$) minus total exports back to the grid ($exp_{c,h}^{P}$) at a certain price for buying ($PBE_{c,h}^{P}$) and another for selling ($PSE_{c,h}^{P}$) energy. If the buying and selling prices are the same, the rate is said to be symmetric in terms of energy. Finally, the parameter $W$ stands for the scaling factor to provide costs on a temporal basis. To calculate the yearly operational costs, $W$ takes the value $\frac{Total\ number\ of\ Hours}{H}$, in which $H$ is the number of hours for the time horizon considered. The net energy costs are then described as:

$$Cost^{P} = \sum_{c \in C} \sum_{h \in H} W_h * \left( imp_{c,h}^{P} * PBE_{c,h}^{P} - exp_{c,h}^{P} * PSE_{c,h}^{P} \right) \tag{3}$$

The last term in the objective function accounts for the network costs, calculated according to Eq. (4). The cost of operating and reinforcing the distribution grid is borne by the distribution system operator (DSO), which is assumed to recover its costs via a regulated grid tariff. We assume two types of costs: sunk costs and prospective costs. The sunk costs are those incurred in the past to build and reinforce the grid in order to meet future demand. Prospective costs are variable and depend on the maximum cumulative load of all customers. The variable *agc* is the additional grid capacity needed to ensure proper functioning with an annualised incremental cost for grid capacity $I^{DSO}$. This can be formulated as in the following equation:

$$Cost^{N} = SunkCosts + I^{DSO} * agc \quad (4)$$

*3.2.2. Grid capacity constraints*

The grid capacity needed to ensure power delivery to customers depends on their hourly imports and exports. Because different customers can import and export energy at the same time, the aggregated net electricity flow is the variable needed to decide the size of additional grid investments. This flow is calculated as the absolute value of aggregate trading by all customers as shown in Eq. (5).

$$e_h^G = \left| \sum_{c \in C} imp_{c,h}^P - exp_{c,h}^P \right|, \forall h \in H \quad (5)$$

Then the total grid capacity, defined as the sum of the existing (*egc*) and additional capacity (*agc*), should always be greater than or equal to the hourly imports (or exports) of all the customers' demands, as formulated in Eq. (6):

$$egc + agc \geqslant e_h^G \; \forall h \in H \quad (6)$$

*3.2.3. DSO cost-recovery constraint*

An important constraint on the upper level is the recovery of grid costs by the DSO via the network tariff. In our framework, the regulator is fully in charge of setting the tariff and we assume that all costs can be recovered via a three-part tariff formulation including volumetric (*vnt*), capacity (*cnt*), and fixed (*fnt*) elements, as in Eq. (7). The three terms respectively depend on each agent's net energy (€/kWh), each consumer's maximum peak (€/kW) and the number of customers *n* (€/customer). Net-metering options are also considered by parameter (*NM*), which can assume different values according to the tariff structure. For example, a user can be charged for the net amount of energy consumed (*NM* = 1), only for imports (*NM* = 0), or for both imports and exports back to the grid (*NM* = −1). This formulation allows the regulator to adopt the tariff which results in the highest social welfare while taking into account customer reactions. The full cost-recovery constraint is described as:

$$Cost^{N} = \sum_{c \in C} \sum_{h \in H} W * \left( imp_{c,h}^P - NM * exp_{c,h}^P \right) * vnt + \sum_{c \in C} p_c * cnt + n * fnt \quad (7)$$

*3.3. Lower-level problem formulation: the agents*

*3.3.1. Agents' objective function*

The objective function of the agents in the lower level is to minimise their total costs subject to the tariff applied by the regulator. It is formulated in Eq. (8) as the sum of investment costs in DER ($Cost^{DER}_c$), energy charges ($Cost^P_c$) and network charges ($Cost^N_c$). The term $Cost^{DER}_c$ is only present for prosumers, as they are able to invest in solar-PV panels and batteries to reduce their peak demand and general energy consumption from the grid.

$$Min\ Cost_c = Cost_c^{DER} + Cost_c^{P} + Cost_c^{N} \tag{8}$$

The following Eqs. (9) to (11) detail each term in the lower-level objective function. First, DER investment costs are the sum of annualised solar-PV and battery investments as in Eq. (9). The energy charges for each customer are related to the type of energy rate that they have adopted from among the retailer's offer, as indicated in Eq. (10). Then, the network charges calculated in Eq. (11) are the charges paid by each customer for grid utilisation. If we consider a dedicated measurement of power and energy for EVs, it will create another connection point which can be physical, using an extra meter, or virtual, via submetering. This is crucial for network tariff design purposes since EVs can then be considered as additional agents separate from the dwelling, even though the dwelling will still bear the costs of charging or collect the revenue from discharging into the grid.

$$Cost_c^{DER} = I_c^S * ic_c^S + I_c^{PV} * ic_c^{PV}, \forall c \in C \tag{9}$$

$$Cost_c^{P} = \sum_{h \in H} W * \left(imp_{c,h}^P * PBE_{c,h}^P - exp_{c,h}^P * PSE_{c,h}^P\right), \forall c \in C \tag{10}$$

$$Cost_c^{N} = \sum_{h \in H} W * \left(imp_{c,h}^P - NM * exp_{c,h}^P\right) * vnt + p_c * cnt + fnt, \forall c \in C \tag{11}$$

*3.3.2. Agents' peak power constraint*

The amount of capacity charges paid by each user depends on their peak power consumption over a time period. As the installed grid capacity must ensure that both bidirectional power flows can be managed, the peak power can occur while customers are importing or exporting energy, as is described in Eq. (12). Since only one term on the left-hand side of the inequality can be non-zero in each time step, we can measure their maximum imported or exported power to calculate the amount of capacity-based charges that should be paid. This is described as:

$$imp_{c,h}^{P} + exp_{c,h}^{P} \leqslant p_c , \forall c \in C, h \in H : \left(\mu_{c,h}^{G}\right) \tag{12}$$

The agents represented by the lower-level optimisation problem are subject to several constraints which are described in Eqs. (13) to (30) and their respective dual variables. These equations describe the investment options available for DERs and their interaction with dwelling load and electric vehicles, and are detailed in the following section.

*3.3.3. Electric vehicle capacity constraints*

The main purpose of EVs is to satisfy the mobility needs of their owners. However, when idle they can be considered as batteries able to inject energy back into the grid while maintaining a certain amount of energy for driving. Eq. (13) describes how the state of charge (SOC) of the EV battery ($s_{c,h}^{EV}$) depends on its state in the previous time step ($s_{c,h-1}^{EV}$), the charging decision ($d_{c,h}^{\Delta EV+}$), the discharging decision ($d_{c,h}^{\Delta EV-}$) and consumption while driving ($D_{c,h}^{\Delta EV-})^2$. Losses in the storage system are represented by the converter loss parameter ($L^S$) and the battery self-discharge parameter ($R_c$). Initial conditions are needed to account for overnight charging and the initial state of charge. In order to enable this, the last time step is linked to the first one, as in Eq. (14), which assumes that the initial and final states of the battery should be equivalent. Finally, the state of charge is determined by parameter $SOC_0^{EV}$ in Eq. (15) below:

$$s_{c,h}^{EV} = S_{c,h-1}^{EV} * \left(1 - R_c^{EV}\right) + d_{c,h}^{\Delta EV+} * \left(1 - L^{EV}\right) - d_{c,h}^{\Delta EV-} * \left(1 + L^{EV}\right) - D_{c,h}^{\Delta EV-}, \forall c \in C, h \in H \setminus \{1\} : \left(\lambda_{c,h}^{EV1}\right) \tag{13}$$

$$s_{c,1}^{EV} = S_{c,H}^{EV} * \left(1 - R_c^{EV}\right) + d_{c,1}^{\Delta EV+} * \left(1 - L^{EV}\right) - d_{c,1}^{\Delta EV-} * \left(1 + L^{EV}\right) - D_{c,1}^{\Delta EV-}, \forall c \in C : \left(\lambda_{c,1}^{EV1}\right) \tag{14}$$

$$s_{c,H}^{EV} = SOC_0^{EV}, \ \forall c, h \quad \left(\lambda_c^{EV2}\right). \tag{15}$$

Eqs. (16) to (19) describe the operational limits of the EV battery. First, Eqs. (16) and (17) ensure that the state of charge of the battery remains within a certain range to avoid extra battery degradation. By limiting the state of charge, we implicitly remove the need to include degradation costs directly in the objective function. Regarding power levels, parameters ($P_{c,h}^{EVch}$) in Eq. (18) for charging and ($P_{c,h}^{EV\,dis}$) in Eq. (19) for discharging are responsible for limiting power transfers, depending on the type of electric vehicle charging equipment (EVSE) adopted. These constraints are described as:

$$s_{c,h}^{EV} \leqslant \overline{E}_{c,h}^{EV},\ \forall c \in C, h \in H : \left(\mu_{c,h}^{EV2}\right) \quad (16)$$

$$s_{c,h}^{EV} \geqslant \underline{E}_{c,h}^{EV},\ \forall c \in C, h \in H : \left(\mu_{c,h}^{EV3}\right) \quad (17)$$

$$d_{c,h}^{\Delta EV+} \leqslant P_{c,h}^{EVch},\ \forall c \in C, h \in H : \left(\mu_{c,h}^{EV4}\right) \quad (18)$$

$$d_{c,h}^{\Delta EV-} \leqslant P_{c,h}^{EV\,dis},\ \forall c \in C, h \in H : \left(\mu_{c,h}^{EV5}\right) \quad (19)$$

[2] This formulation can potentially lead to discharging and charging episodes happening in the same time step. The solutions are verified ex-post to ensure that realistic behaviour happens for all agents over the entire time horizon.

*3.3.4. Battery storage investment constraints*

A stationary battery allows the user to temporarily shift load and store the surplus electricity generated by solar PV. Analogous to the EV charging and discharging equations, Eq. (20) describes how the state of charge of the battery ($s_{c,h}$) depends on its state in the previous time step ($s_{c,h-1}$), the charging decision ($d^{\Delta+}_{c,h}$) and the discharging decision ($d^{\Delta}_{c,h}{}^{-}$). In this case, we let the optimisation define the initial state of charge of the battery since it will depend on the total size of the battery installed, which is also a decision variable.

$$s_{c,h} = s_{c,h-1} * (1 - R_c) + d_{c,h}^{\Delta+} * \left(1 - L^S\right) - d_{c,h}^{\Delta-} * \left(1 + L^S\right),\ \forall c \in C, h \in H \setminus \{1\} : \left(\lambda_{c,h}^{S1}\right) \quad (20)$$

$$s_{c,1} = s_{c,H} * (1 - R_c) + d_1^{\Delta+} * \left(1 - L^S\right) - d_1^{\Delta-} * \left(1 + L^S\right),\ \forall c \in C : \left(\lambda_{c,1}^{S1}\right) \quad (21)$$

The capacity installed is decided in the model, if it is profitable for each agent, by choosing variable $ic_c^s$ bounded by a maximum capacity limit ($U_c^s$), as in Eq. (22). For certain agents who do not have the possibility of installing any battery capacity, the maximum value can be set at zero.

$$ic_c^s \leqslant U_c^S, \forall c \in C : (\mu_c^{S1}) \tag{22}$$

As in the EV case, Eqs. (23) to (26) describe the operational limits for the stationary battery. The parameters ($S^{\%max}$) and ($S^{\%min}$) are the percentages of maximum and minimum charge levels allowed respectively. Regarding power levels, the charging factor ($P_c^{ch}$) in Eq. (25) and discharging factor ($P_c^{dis}$) in Eq. (26) represent the maximum limits of power transfer according to the storage system specifications, described as follows:

(23)

$$s_{c,h} \geqslant S^{\%min} * ic_c^s, \ \forall c \in C, h \in H : (\mu_{c,h}^{S3})$$

$$d_{c,h}^{\Delta+} \leqslant ic_c^s * P_c^{ch}, \ \forall c \in C, h \in H : (\mu_{c,h}^{S4})$$

$$d_{c,h}^{\Delta-} \leqslant ic_c^s * P_c^{dis}, \ \forall c \in C, h \in H : (\mu_{c,h}^{S5})$$

(24)

(25)

(26)

*3.3.5. Solar-PV investment constraints*

The solar-PV capacity installed is also endogenously decided in the model. If it is profitable for each agent, this is done through variable $ic_c^{PV}$, bounded by a maximum capacity limit ($U_c^{PV}$) as in Eq. (27). For certain agents incapable of installing any solar PV, the maximum value can be set at zero. The amount of energy produced will depend on the solar availability in kW/kWp ($G^{PV}_{c,h}$) and users do not have the option of curtailing, meaning that they have to export the surplus electricity generated in any given time period.

$$ic_c^{PV} \leqslant U_c^{PV} \forall c \in C : (\mu_{c,h}^{PV1}) \tag{27}$$

*3.3.6. Energy balance equation*

The energy balance equality couples all the investment and operational decisions with the load demand profile of each customer ($D_{c,h}$), as in Eq. (28). The terms $imp^L_{c,h}$ and $exp^L_{c,h}$ allow the interaction between an EV which is metered separately and the household load to be modelled. By treating the EV as a separate agent, it can buy electricity at an EV-only rate and if needed transfer power to the house, which is subject to a different rate. In this manner, Eq. (29) ensures the supply-demand balance: all the imports are equal to the exports. Similarly, households that may invest in DERs will have their own interaction with their electric vehicles. Now, dwellings can invest in DERs and arbitrage energy between the two different retail tariffs, for example by charging the EV with a battery at a flat rate. It is also possible to charge the EV using local solar energy instead of relying on buying electricity from the grid, significantly increasing the complexity of the interactions. To restrain the relation between EVs and prosumer-type households, another equilibrium equation (Eq. (30)) is added. The auxiliary parameter $\alpha_c$ forbids other agent types from interacting with their dwellings at EV-only rates by setting it to zero. An analogy can be made with local market modelling. In our case, an EV and a specific type of household form a local market in which they can only interact with others via the main grid. These interactions are described in Eqs. (28) to (30):

$$D_{c,h} + d_{c,h}^{\Delta EV+} - d_{c,h}^{\Delta EV-} + d_{c,h}^{\Delta +} - d_{c,h}^{\Delta -} - ic_c^{PV} * G_{c,h}^{PV} = imp_{c,h}^{P} - exp_{c,h}^{P} + \alpha_c * \left(imp_{c,h}^{L} - exp_{c,h}^{L}\right),\ \forall c \in C, h \in H\ :\ \left(\lambda_{c,h}^{EB}\right) \quad (28)$$

$$\sum_{c \in C_{EV}} \left(imp_{c,h}^{L} - exp_{c,h}^{L}\right) = 0\ , \forall h \in H\ :\ \left(\lambda_h^{L_{EV}}\right) \quad (29)$$

$$\sum_{c \in C_{EV/DER}} \left(imp_{c,h}^{L} - exp_{c,h}^{L}\right) = 0\ , \forall h \in H\ :\ \left(\lambda_h^{L_{EV/DER}}\right) \quad (30)$$

Additional constraints can be added to enforce the interaction between EVs and the houses to which they are connected, as in Eqs. (31) and (32). These equations limit the amount of energy that a vehicle can import from the house according to the energy resources existing (battery and PV). Reciprocally, the house can also only import energy from the vehicle if it is ready to discharge at any time. The model allows V2H to be avoided by forcing the constraint of equation Eq. (32) to zero while allowing V2G via the energy balance equation (28).

$$imp_{EV,h}^{L} \leqslant d_{c,h}^{\Delta -} + ic_c^{PV} * G_{c,h}^{PV},\ \forall c \in (C_{EV} \cup C_{EV/DER})\ ,\ \forall h \in H : \left(\mu_{c,h}^{impL}\right) \quad (31)$$

$$imp_{c,h}^{L} \leqslant d_{EV,h}^{\Delta EV-},\ \forall c \in C_{EV}, \forall c \in (C_{EV} \cup C_{EV/DER})\ ,\ \forall h \in H\ : \left(\mu_{c,h}^{impL2}\right) \quad (32)$$

## 4. Case study: setup and input data

In this section, the setup and input data for a case study using the bi-level model are described. First, the general setup of the numerical example will be explained as a starting point. Then, the data regarding agents such as load profiles, solar insolation and energy tariff profiles are presented alongside the data on EVs and DERs. Finally, the baseline electricity bill is followed by the grid cost structure.

*4.1. Setup*

In the case study, two behaviours are considered regarding the option that agents have to invest in DER like solar PV and stationary batteries: prosumer and consumer behaviour. Moreover, when EVs are considered, a combination resulting in four different types of agents is observed: prosumer with EV, consumer, prosumer and consumer with EV.[3] Given that smart charging, V2G and submetering are considered, the number of possible scenarios could rapidly increase. To limit the number of agents and scenarios, V2G is only adopted by prosumers with EVs since they can be considered full innovators and less risk-averse. Submetering is adopted for those who can obtain higher gains according to the first MILP model results.[4] Smart charging is essential to limit the risk of surpassing the maximum capacity of the grid during the peak. Therefore, it is considered a common practice among all EV owners. The representation of the agents depends on which type of EV rate they are charged at and if they are able to make DER investments, as is shown in Fig. 2:

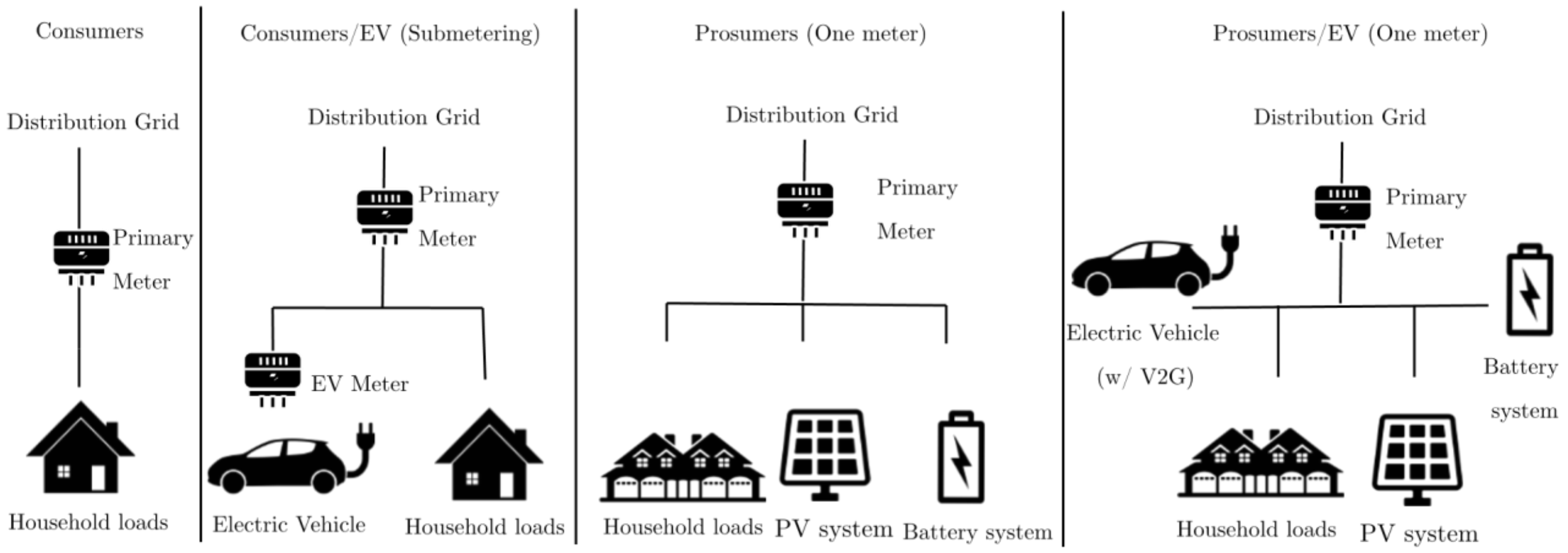

Figure 2: Topology configuration of all agents.

The topology configuration of agents complexifies when moving towards prosumers with submetered EVs. In this case, one rate is specifically applied for EV charging using the existing EVSE meter and another for the remaining loads, including DERs. It is possible, in theory, to arbitrage energy within the same household, collecting the spread between the rates on condition that investments in DERs or V2G are made. The assessment using the MILP model considers different energy profiles, topologies and DER costs, so EV owners could choose the options that reduce their total cost.

[3] In this configuration a 50% share of prosumers is analysed. Although this seems a high share at today's global level, in a mid-term perspective the number of prosumer agents is expected to increase. Increasing environmental awareness and DER technology cost reductions contribute to boosting the number of consumers becoming prosumers.

[4] With reference to the share of EVs, 50% of agents are considered EV owners. This proportion is relative to the total vehicle stock present if the remaining agents are considered to have ICE vehicles. The same factors influencing the shift to prosumer behaviour will have an important impact on electric mobility uptake. Moreover, current policy support for zero-emission light-duty vehicles and ICE car bans in more than 20 countries by 2030-2035 will also contribute to a high penetration of EVs in the vehicle market (IEA, 2021).

The first analysis will consider only different energy profiles (flat or time-of-use) while the network charges are fixed according to the baseline electricity bill. This baseline scenario will be used as a reference for all the other counterfactual scenarios with respect to energy and network cost variation. The grid costs expected to be recovered via tariffs depend on the state of the grid. In other words, if all the investments have already been made the costs are considered to be sunk, or if there are still investments to be made they are considered to be prospective. Once the model is calibrated with input data, the MPEC formulation will allow evaluation of the variation in grid charges for agents according to how the regulator sets the tariffs.

*4.2. Load and solar profiles*

The 48-hour load profiles adopted correspond to the inelastic hourly demand by prosumers who have the capacity and the means to invest in DERs and by consumers who are not able to. A total annual electricity consumption of 10,000 kWh is adopted for prosumers, which is equivalent to the average consumption of a U.S. residential utility customer (EIA, 2020). A value of 5,500 kWh is chosen for consumers, being equivalent to average French dwelling consumption (Odyssée Mure, 2019). This difference in total consumption is explained by a positive correlation frequently observed in the literature between income and electricity use (Burger et al., 2020; Borenstein, 2012). Peak demand values for a typical household tend to occur in the early evening with a second smaller peak in the morning. Here, they differ for the two types of agents, a peak of 4.8 kW for prosumers and 3.2 kW for consumers, which leads to a coincidental peak of 8 kW without any EV or DER investment.

Regarding the solar profile, a 48-hour profile with two different insolation peaks is used. The first one with higher insolation represents a typical sunny day in which there is direct sunlight without any external interference in all periods. The second represents a cloudy or rainy day when the solar irradiation is deeply reduced, leading to less solar PV electricity production. A day with a smaller peak of insolation is synchronised with the days presenting higher peaks of electricity consumption. Both the load and solar profiles are illustrated in Fig. 3.

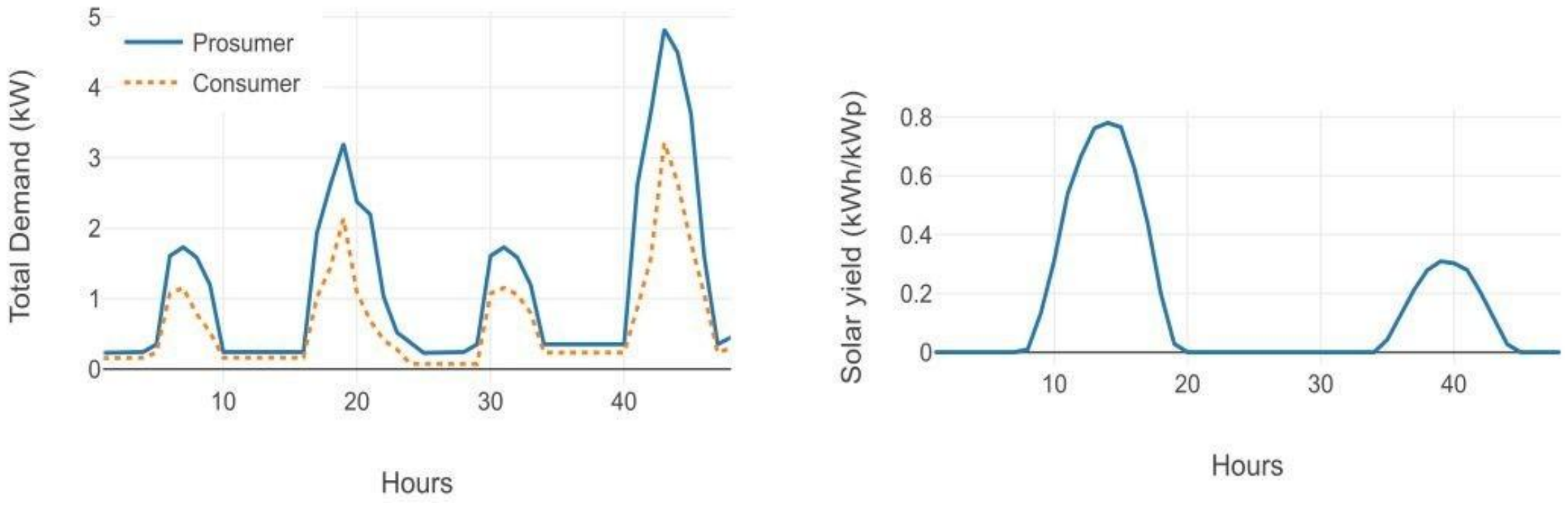

Figure 3: Load and solar profiles.

*4.3. EV and DER data*

Electric vehicle charging constitutes a large share of the residential electricity bill, especially for owners with important mobility needs. According to a pilot conducted in the United Kingdom exploring user behaviour related to home charging, 3,500 kWh corresponds to the average yearly consumption of vehicles with batteries greater than 35 kWh (Western Power Distribution, 2019). This amount is equivalent to the average electricity consumption of a dwelling in the U.K. An EVSE allowing a maximum input power of 7 kW is adopted. Then the peak caused by EVs will depend on the type of rate applied since it is possible to spread the charging over the day using smart charging.

From the yearly EV energy consumption it is possible to derive daily energy needs for mobility purposes. First, based on an average energy consumption of 194 Wh/km (Electric Vehicle Database, 2021a), a daily distance travelled of around 50 km/day is calculated.[5] This distance is in line with that found for EVs in European countries like Italy, France and Germany (Pasaoglu et al., 2014) and in California (California Energy Commission, 2019). Therefore, a daily energy need of 9.7 kWh is estimated for the case study based on these data.[6]

To avoid extra battery degradation, it is advisable to keep the battery state of charge in a specific range and to have a limited depth of discharge. It is assumed that the battery capacity of fully electric vehicles has an average value of 60 kWh (Electric Vehicle Database, 2021b) and the SOC is allowed to vary between 6 and 54 kWh (10-90% respectively). Finally, an assumption regarding connection hours is defined as observed in the British EV charging pilot. The majority of EVs disconnect from their homes at 7am and plug back in at 5pm (see Appendix B for more details). Conversion losses from electronic power converters in the EVSE of 5% and a negligible self-discharge rate complete the EV technical parameters.

The cost of technologies such as solar PV has a strong impact on the size of agents' DER investments and vary substantially across the utility, commercial, and residential sectors (NREL, 2020). The cost adopted for installing solar PV is around 900 $/kWp, which can be considered low in the context of residential-scale systems. However, we use this optimistic scenario to illustrate a situation where it is cost-optimal for agents to invest in these systems. Direct subsidies can also support the argument for using this cost for end-users. On the other hand, net metering will not be considered since it may over-incentivise solar PV investments and cause strong fairness issues among agents. A discount factor of 5% and a lifetime of 20 years translates into an annualised cost of 72 $/kWp.

Stationary batteries have been benefiting from the decrease in the cost of automotive batteries at the pack level. The cost adopted for them is 150 $/kWh, which corresponds to the current weighted average cost of a

---

[5] $Daily\ \ distance\ travelled = \frac{Yearly\ average\ consumption}{Efficiency * Days\ in\ a\ year} = \frac{3{,}500}{194\frac{Wh}{km} * 360\ days} = \ 50.1\frac{km}{day}$

[6] $Daily\ energy\ needs\ = \frac{Yearly\ average\ consumption}{Days\ in\ a\ year} = \frac{3{,}500}{360\ days} = \ 9.7\frac{kWh}{day}$

battery pack for electric vehicles (IEA, 2021). In addition, an annualised cost of 19.4 $/kWh is obtained by using a similar discount factor of 5% and a lifetime of 10 years. Regarding the technical parameters, we assume conversion losses of 5% and a self-discharge rate of 0.1% per hour.[7] Finally, analogous to the EV battery case, the SOC is allowed to vary between 10% and 90% of total battery capacity.

### *4.4. Baseline electricity bill*

A reference scenario is needed to compare the outcome of the optimisation model for further analysis. The baseline electricity bill is a counterfactual bill defining energy and network costs in dollars for each agent. Furthermore, renewable energy support (RES) and taxes are also important components in the majority of electricity offers for households. The weighted average breakdown of electricity offers in 2019 provided by the ACER marketing monitoring report (ACER, 2020) indicates a 45% share of energy costs and 33% of network charges. If the amount of RES is split equally between them, taxes will account for the remaining 22%. Electricity prices vary considerably across countries in different regions. For instance, the average price in the E.U. is 21.6 euro cents/kWh (ACER, 2020) while in the U.S. the average was found to be 10.6 cents/kWh (EIA, 2019). An average value of 16.3 cents/kWh is used as the reference for the bill calculation coupled with the breakdown information. Table 1 shows the final bill for each agent in detail.

Table 1: Electricity bill components

| | Breakdown | Cost in bill ($/kWh) | Consumer ($) | Consumer/EV ($) | Prosumer ($) | Prosumer/EV ($) |
|---|---|---|---|---|---|---|
| Energy | 45% | 0.073035 | 406 | 662 | 731 | 986 |
| Network | 33% | 0.053559 | 298 | 485 | 536 | 723 |
| Other charges | 22% | 0.035706 | 198 | 323 | 357 | 482 |
| Total costs | 100% | 0.163 | 902 | 1470 | 1624 | 2192 |

So far, energy costs have been assumed to be invariant over time, meaning that a flat energy profile is adopted. Nevertheless, this type of energy profile is not suitable to incentivise households to reduce their consumption during certain periods of the day or to install DERs for peak-shaving and valley-filling purposes. Time-of-use energy profiles can be the appropriate economic signal to meet these objectives. Therefore, two extra TOU profiles are proposed for evaluation. The first, denominated TOU1, has the highest relative value synchronised with the private peak, incentivising agents to offset coincidental demand and reduce their own peaks. The second profile, TOU2, supports solar PV adoption since its relative value around the period of solar production is the greatest among all profiles. In order to make them comparable, both TOU profiles are calibrated so that the final consumer bill in the base case scenario is the same regardless of the choice of energy profiles. The final average price of TOU rates is often lower than the equivalent flat rate due to the

[7] The presence of self-discharge only in stationary battery systems and not for EV batteries results from an assumption that battery management systems (BMS) in EVs are more efficient. By not allowing the battery to function outside its operating margins with respect to temperature, for example, self-discharge will be greatly reduced, therefore remaining negligible for functioning purposes.

high number of mid-peak and off-peak hours. However, the final electricity bill is slightly higher. This shows how complex tuning energy profiles can be, depending mostly on what kind of incentives the utility purposes. For instance, it may increase on-peak charges to avoid grid congestion or decrease the off-peak level to encourage users to shift their consumption. We note that the rates selected are not symmetric, which means that the compensation for injecting electricity into the grid is not at the same level as the cost of withdrawing electricity. Users are subject to 10% compensation related to the energy price at the moment of injection according to the type of tariff.[8] All the energy profiles are illustrated in Fig. 4:

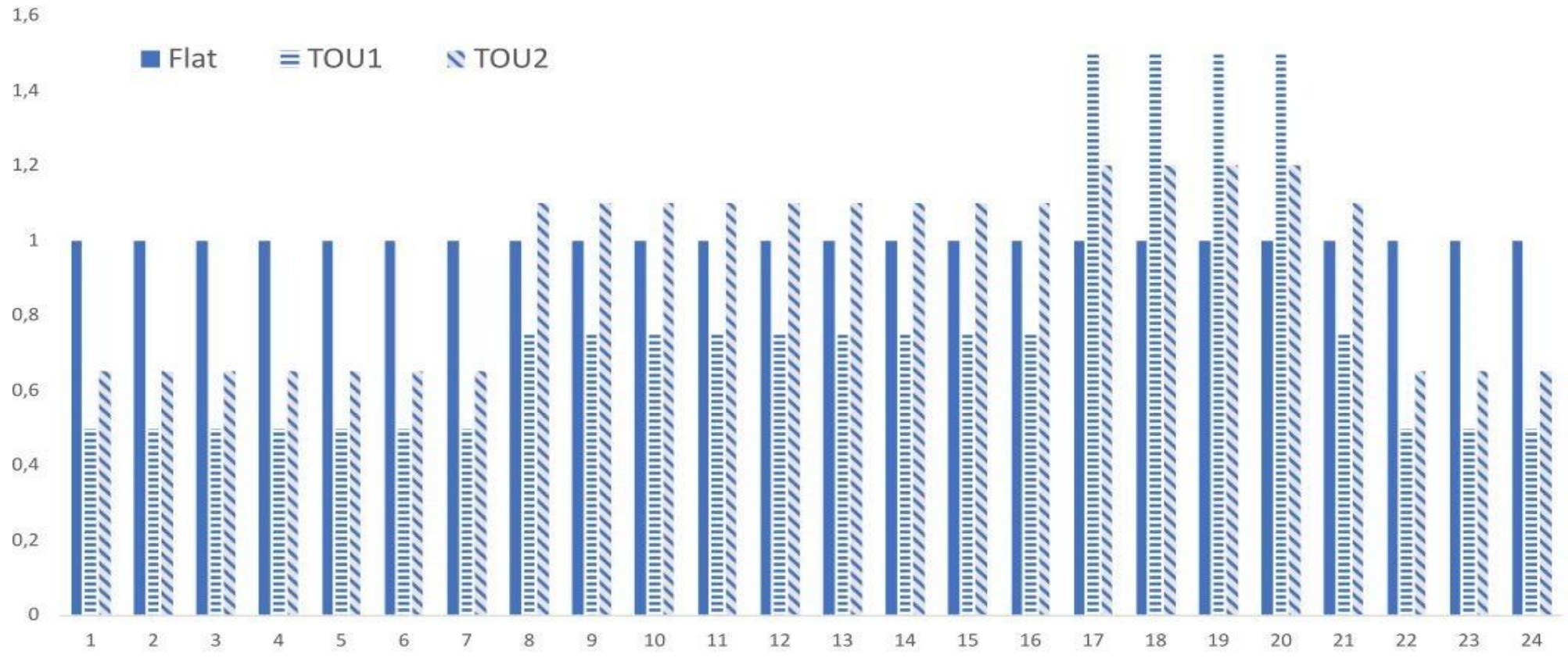

Figure 4: Energy profiles.

*4.5. Grid cost structure*

Several factors affect the grid cost structure, which makes assumptions regarding the structure complex to make. For instance, fixed operating costs, sunk capital costs, variable operating costs and network losses are all parts of the cost needed to be recovered by network operators. According to Simshauser and Downer (2016), in an electricity distribution system the fixed and sunk capital costs will comprise 70-80% of the total cost structure. This indicates that traditional distribution grids are over-dimensioned for current demand, a practice which is described as 'fit-and-forget.'

Nevertheless, the penetration of new distributed energy resources (stationary batteries, solar PV and heat pumps) and electric vehicles will substantially change this scenario. Grid reinforcements may be necessary to cope with the peak load increase caused by these technologies if their use is not properly coordinated. The economic feasibility of substituting grid capacity with local flexibility can be assessed by using forward-looking grid costs. A reduction of network utilisation cannot reduce the costs of the current network infrastructure, which are sunk, but only defer future network investments by reducing coincident peak loads

[8] The future of policies regarding PV injection prices is quite uncertain. Due to high cost-shifting levels in some locations, indirect incentives to instal PVs such as feed-in tariffs and compensation for energy injection are thought to substantially decrease (IOUs, 2021).

(Govaerts et al., 2019). Therefore, this element, also named prospective costs, represents the long-run costs of the network for agents.

As in Schittekatte and Meeus (2020), we adopt three scenarios for the analysis: 100% sunk costs; 50% sunk and 50% prospective; and 100% prospective grid costs. For the first scenario, the *SunkCosts* parameter in Eq. (4) is set at $ 2,040, which is the sum of all agents' network charges in the baseline electricity bill (see Table 1). For the second scenario, half of these costs are employed as sunk costs and the annualised incremental cost for grid capacity, $I^{DSO}$, representing the prospective costs is 63.8 $/kW.[9] Finally, in the full prospective scenario the prospective cost is set at the annualised incremental cost for grid capacity of 127 $/kW while there are no sunk costs.[10] Even though grid costs are lumpy and vary depending on site-specific properties, these values allow us to have a fair representation of grid cost structure.

## 5. Results

In this section we present the results of the numerical model. We begin by analysing the impact of different EV owner energy profiles on the final electricity bill with a fixed base case volumetric network tariff. Then, once the energy profiles are set for the agents owning an EV, we assess the effects of different network tariff designs on the final bill. We end with a discussion and by highlighting some policy implications.

### *5.1. Energy profile assessment*

A retailer's electricity offer usually has several types of rates given the difference in time granularity for customers. Flat or time-varying rates are proposed for residential customers, so they are encouraged to compare current energy costs to other tariff rate options and select the best rate plan (SCE, 2019). The aim of this section is to retain the cost-efficient solution, i.e. one meter or submetering, by choosing a combination of energy profiles, while grid charges are exogenously fixed beforehand to minimise energy costs. The logic behind this choice is to anticipate EV owner actions regarding the choice of rate, assuming that a single consumer decision does not have an impact on other consumers' choices, and therefore has a negligible impact on the network tariff design.

Agents possessing an EV are analysed individually using the MILP model described in Section 3.1 to ascertain their total costs. All three energy profiles are adopted for agents using only one meter (flat, TOU1, TOU2), while for the submetering configuration a combination of a flat rate for house loads and TOU1 or TOU2 for EV charging are considered. According to our setup, only prosumers have the option of adopting V2G technologies since they are considered full innovators. We separate the analysis of consumers who are not able to invest in DER from that of prosumers who can. Figs. 5 and 6 show the variation in the total cost for consumers and prosumers owning an EV under different energy profiles.

[9] $I^{DSO}\,(50\%\,prospective) = \frac{0.5 * Average\ sunk\ costs}{Average\ coincident\ peak} = \frac{0.5 * 510}{4\ kW} = 63.7\,\frac{\$}{kW}$

[10] We assume that EV charging done via smart charging will not deeply affect the coincidental peak demand. Therefore, the average peak used to calculate prospective costs corresponds to the coincident peak without any EV or DER installed.

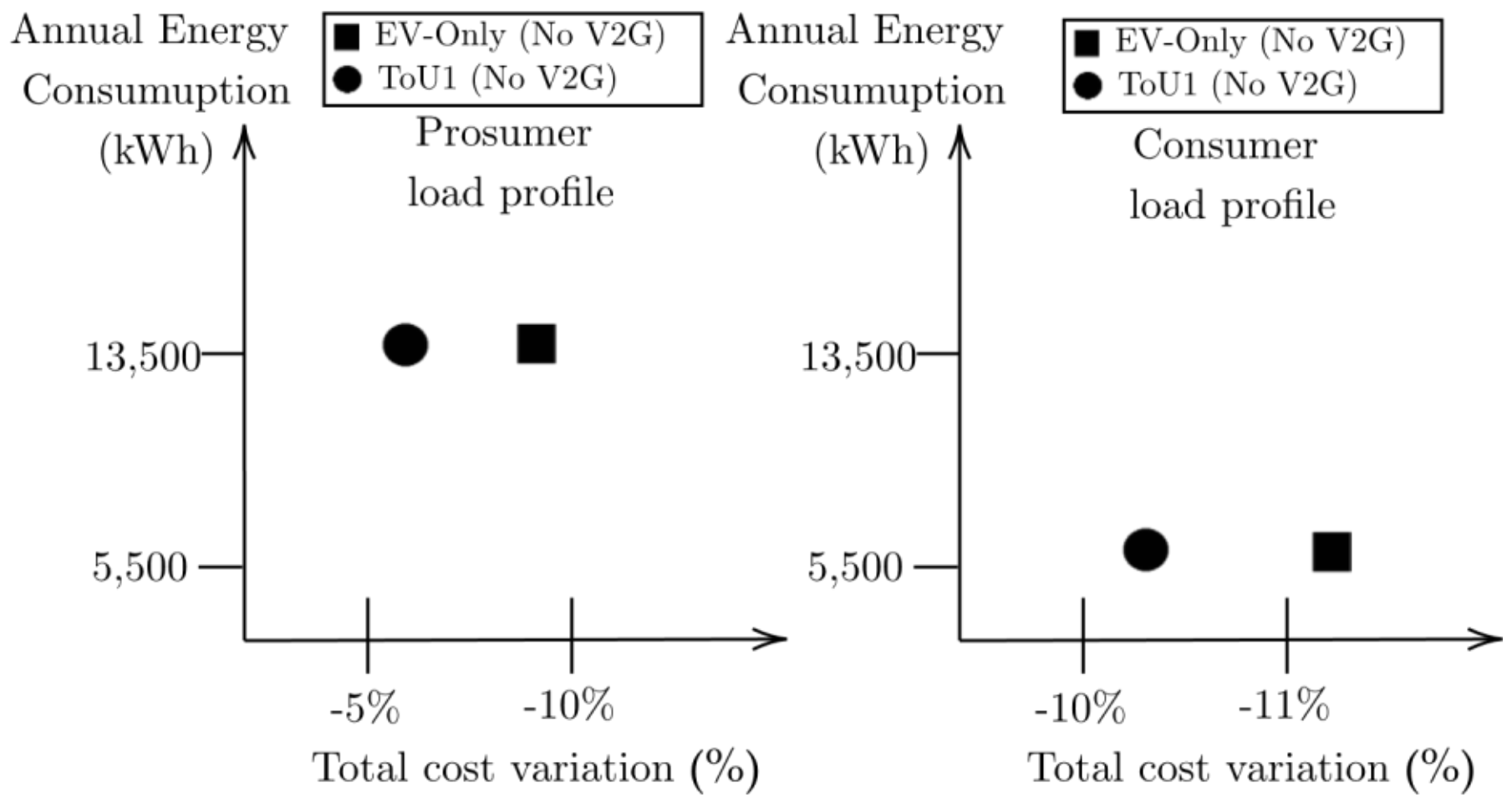

Figure 5: Energy cost variation without DER investment.

First, Fig. 5 shows the cost variation for two rate choices compared to the flat baseline rate for two different load profiles. This clearly reveals the positive impact of submetering by applying EV-only tariffs for agents who are not able to invest in DER. For the prosumer load profile, the cost reduction obtained by adopting an EV-only solution (Flat rate for the house plus TOU1 for the EV) surpasses the whole-house TOU1 tariff by around 2%, while for the consumer profile half of this value is obtained. In conclusion, the greater the energy consumption of the household, the higher are the bill savings brought by the submetering solution.

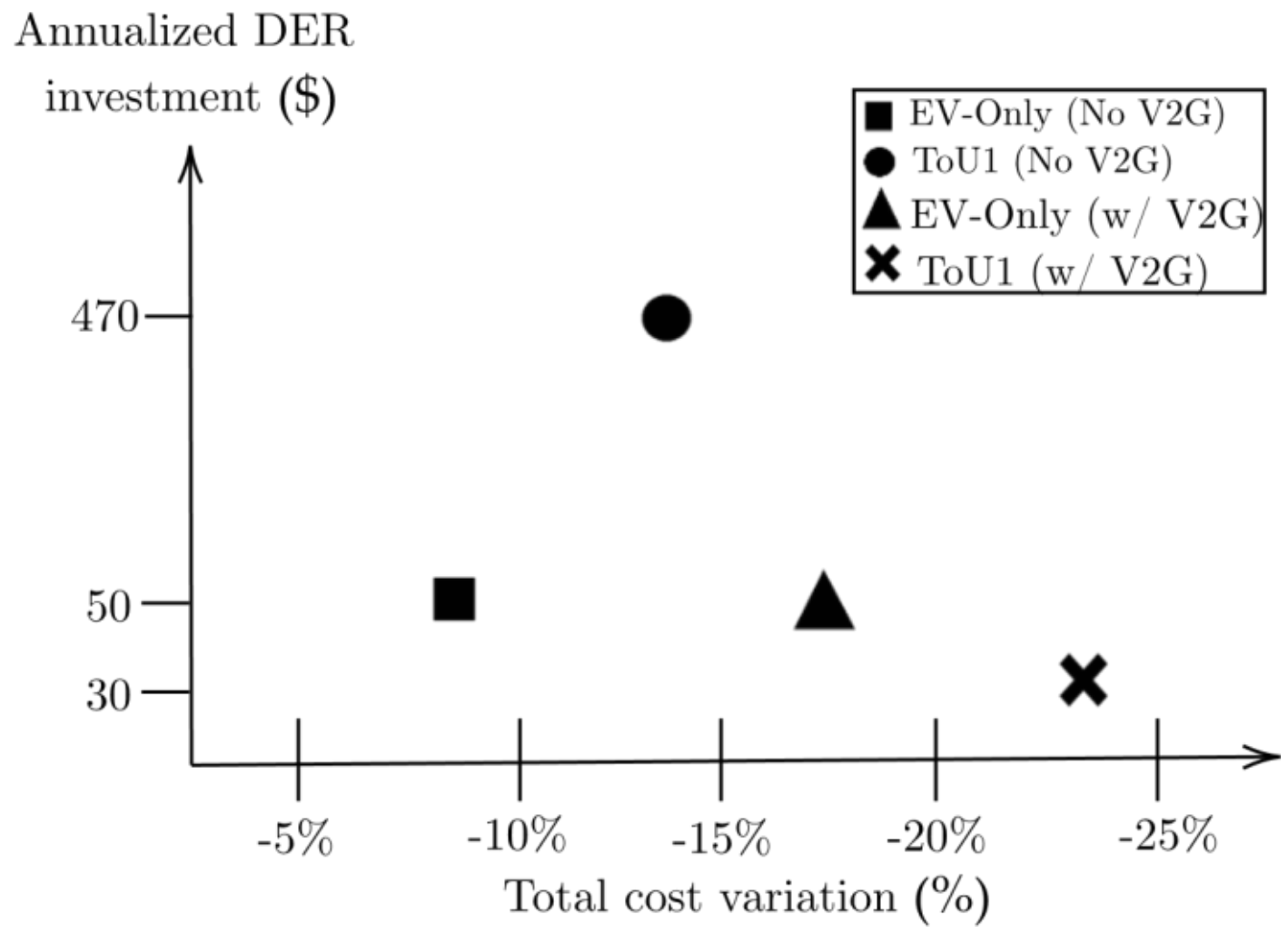

Figure 6: Energy cost variation with DER investment.

Prosumers who are ready to invest in DER may have different choices regarding the type of solution to adopt. Fig. 6 shows that whole-house TOU tariffs are preferred since they bring a greater energy cost reduction with or without V2G. This preference is due to the higher cost-effectiveness of arbitraging energy from off-peak to peak hours using a stationary battery or V2G under the whole-house TOU1. On the other hand, EV-only tariffs do not incentivise the adoption of DERs, especially batteries, considering that they will be charged at a flat rate which is not financially attractive. A great share of spread is already obtained by adopting the different tariffs for the submetering solution, which reduces the gains brought by DERs.[11] The detailed results for both solutions are shown in Table C.2.

Finally, based on these results we are able to choose the energy profiles and type of solution for each agent who owns an EV for further analysis. For prosumers with EVs and V2G, the whole-house TOU1 tariff is adopted as it results in the lowest costs among the options. For consumers who possess an EV, the submetering solution coupling a flat rate for house loads and TOU1 for EV charging is adopted for the same reason of cost-effectiveness. For the comparison scenario, all agents able to invest in DER adopt TOU1 and passive consumers select flat rates.

### *5.2. Network cost impact*

With the energy profiles defined for all agents, it is possible to identify the impacts on network costs caused by EVs and submetering solutions. The MPEC model allows simulation of several tariff designs that can be used by the network operator to recover grid costs. We adopt two possible tariffs structures for this purpose: a pure volumetric tariff, which is the most common way of billing customers in the current residential sector and a three-part tariff including capacity and fixed charges. A technology scenario using V2G technology and submetering is also included. This represents a total of four scenarios which include the two different tariff structures and the possibility of adopting EV submetering coupled with V2G.

First, the changes in total network costs are presented for all three grid cost structure scenarios. Then, individual cost variations and their manner of allocation according to different tariff designs are explored. We focus on the submetering solution including the variation that it may cause in network costs to verify its cost-effectiveness. The cumulative network costs for all scenarios based on three grid cost structures are presented below in Fig. 7.

[11] This scenario could differ according to price. For instance, simulations with a more conservative PV cost of 1,500 $/kWp indicate that EV-only tariffs are more cost-effective than ToU, without a need for great investments in DER. The full data for this scenario are presented in Table C.3.

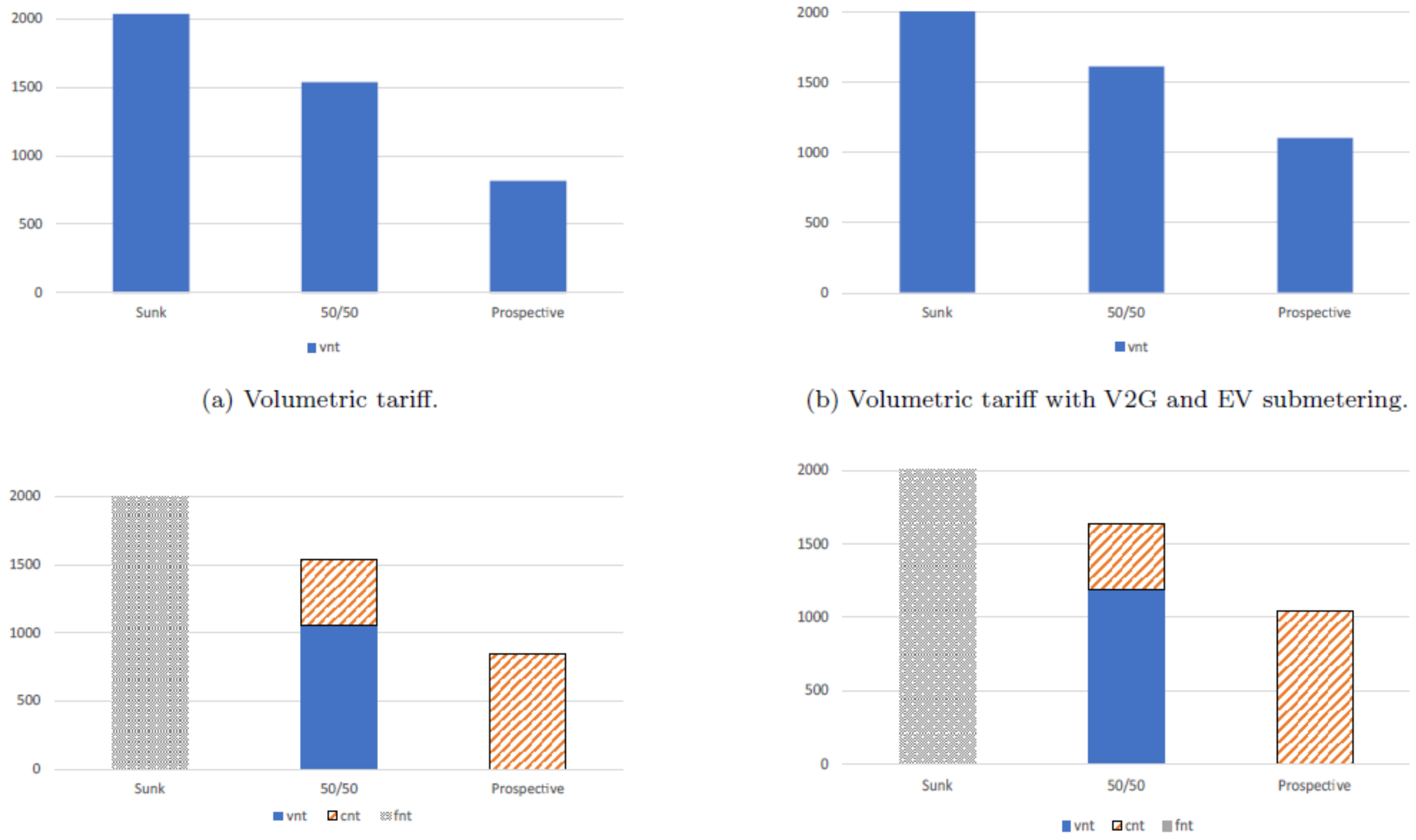

(a) Volumetric tariff.

(b) Volumetric tariff with V2G and EV submetering.

(c) Three-part tariff.

(d) Three-part tariff with V2G and EV submetering.

Figure 7: Total network costs.

Our main observation is that the higher the share of prospective costs in the network cost structure, the lower the total grid costs. This is due to the implicit cooperative behaviour of all agents once variable costs are present and dependent on the coincident peak. In the case of 100% sunk costs, prosumers invest as much as possible in DERs to avoid grid charges and shift more of them to consumers. We focus on the scenarios with prospective costs, given the exhaustive sunk costs analysis in Schittekatte et al. (2018) and Hoarau and Perez (2019). Nevertheless, in both cases with prospective charges prosumers invest in DERs to reduce their contribution to coincident peak increases, while EV owners also use smart charging to avoid constrained hours. As a consequence, the tariff will decrease for all agents since the coincident peak that drives network variable costs is greatly reduced.

We note a perceptible increase in the cumulative network charges for the scenarios with V2G and submetering containing variable charges (see Table C.4 for more numerical details). The underlying logic is that V2G greatly reduces the need for battery investments. Unlike stationary batteries, EVs are away from home during a great part of the day. Therefore, the required charging episodes for driving needs share the same peak period as those required for energy arbitrage. Consequently, an increase in the coincident peak occurring in off-peak periods increases the network tariff. A possible solution to alleviate a coincident peak increase during off-peak periods involves charging EV at work or using public charging infrastructure. If this type of charging is incentivised, users could spread EV charging over all the hours of the day, which would have a positive impact on the total grid costs to be recovered.

The shift from a volumetric to a three-part tariff does not significantly change the network costs, although it impacts users in its allocation. Prospective costs, if recovered via capacity charges, will benefit those who can reduce their private peak by investing in storage or adopting V2G technology. With a volumetric tariff storage will be mostly used to arbitrage energy from on-peak to off-peak periods and lower the coincident peak. This will optimise the energy usage of prosumers without deeply affecting the charges paid by consumers. There will then be fewer fairness issues concerning the volumetric tariff, for which consumers will be spared high cost-shifting levels, as is shown in Fig. 8.

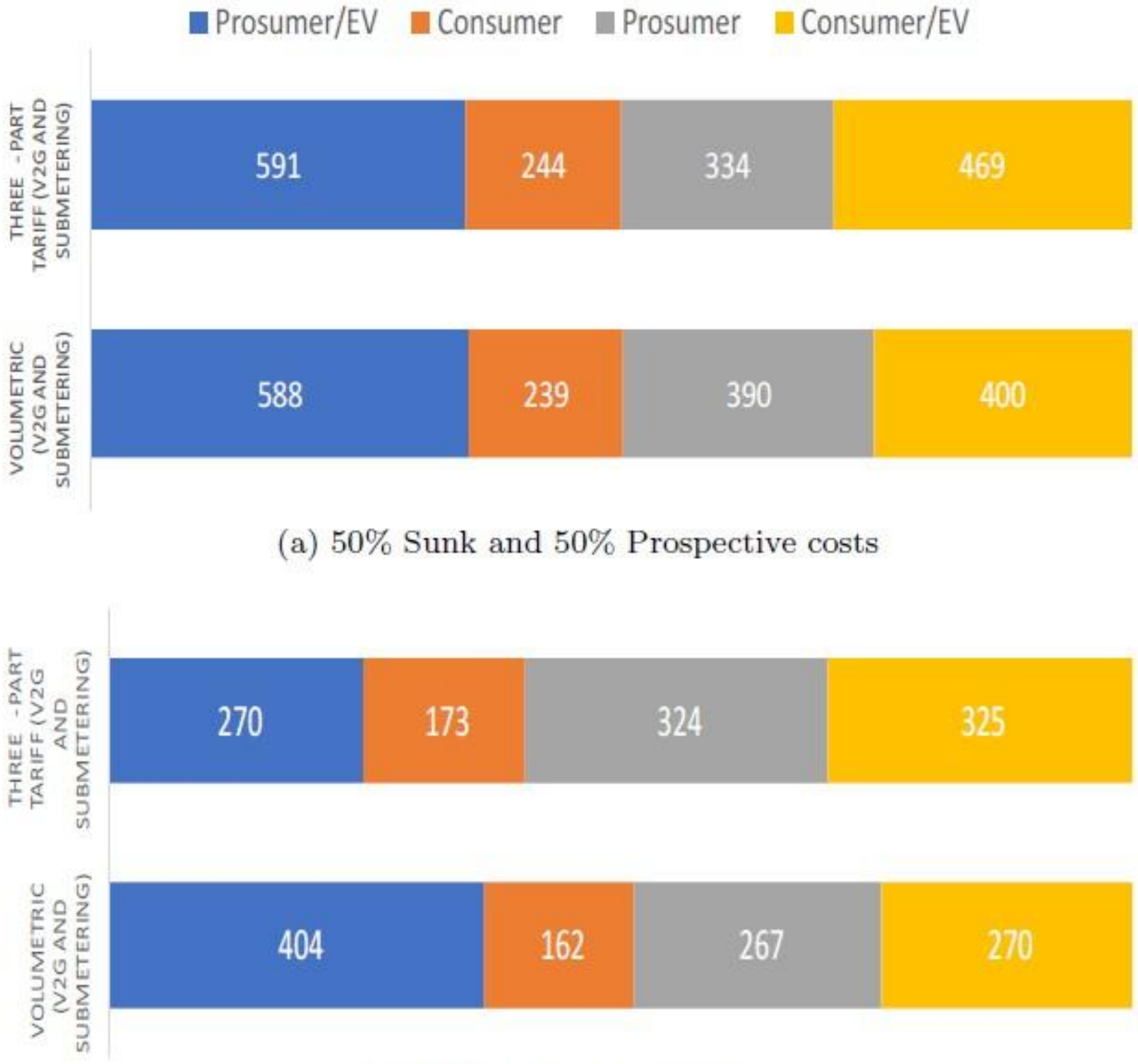

(a) 50% Sunk and 50% Prospective costs

(b) 100% Prospective costs

Figure 8: Network cost breakdown.

Comparing the volumetric tariff with a three-part tariff, consumers experience a 1.5% and 3.7% reduction in grid charges in the 50% sunk/50% prospective and 100% prospective scenarios respectively. The same trend is observed for consumers with an EV, but in this case the magnitude of the reductions are around 14% and 11% respectively. Clearly, a volumetric tariff favours agents possessing an EV (submetered or not) who cannot invest in DERs. Nevertheless, it is essential to determine the total charges for these agents to evaluate the cost-efficiency of the submetering solution. Gains obtained from energy savings could be outweighed by increases in network charges due to other agent's reactions to the type of tariff applied. Therefore, Fig. 9 shows the total costs for consumers with an EV with different tariffs and in different technological scenarios.

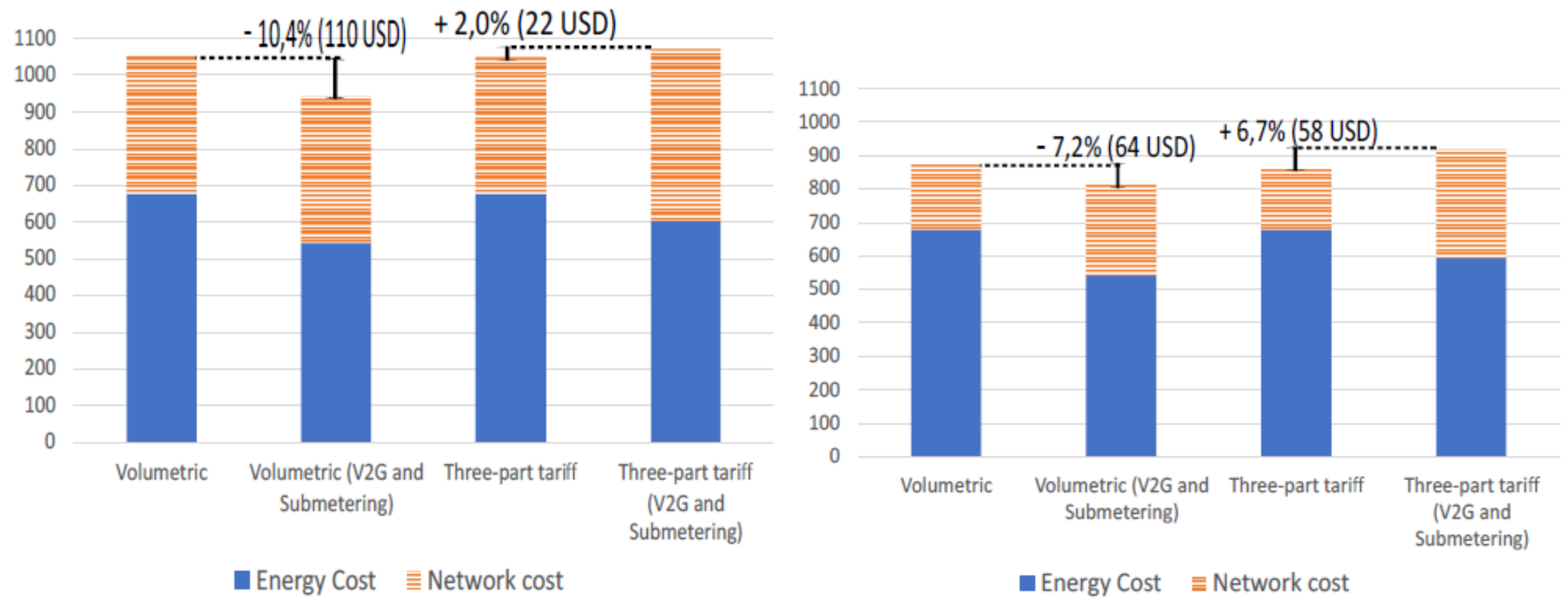

(a) 50% Sunk and 50% Prospective

(b) 100% Prospective

Figure 9: Total costs for consumer/EV agents in all four scenarios.

A dedicated tariff for EVs excluding house loads creates a positive spread of energy charges as was previously described in Section 5.1. An extended analysis which includes the variation in network costs caused by other agent responses in this context will shed light on the potential of this solution. With the volumetric tariff, the rise in network costs driven by a slight increase in the coincident peak does not outrun the reduction in energy costs. The final bill reductions observed are around 10.4% and 7.2% for the 50% sunk/50% prospective and 100% prospective scenarios respectively. In the case of three-part tariffs, the outcome is not so favourable. The capacity charges applied to agents transfer the network savings by prosumers to consumers, including those with an EV with submetering. In this case, the rise in network charges surpasses the energy savings by having two different tariffs. Therefore, final rises of 2.0% and 6.7% are observed.

## 6. Discussion and policy implications

The integrated assessment of submetering exploring energy and network charge variation has shown that this solution brings added value to EV charging. Submetering simplifies the indirect load control of EVs by using a specific economic signal for charging while leaving other appliances charged at another rate. Adoption of an adequate smart meter placed upstream that can properly communicate with the EVSE meter is imperative to achieve this. Pilots and demonstration projects are essential not only to prove the feasibility of the concept but also to establish future protocols and adjust the technology of the appliances involved. The transaction costs involved in adjusting and optimising the use of an EVSE meter would be counterbalanced in a short period, taking into consideration the variety of services potentially available. This configuration of meters also allows other EV energy services, such as frequency regulation, that can be done using V2G without the need to buy or rent a dedicated meter. Stacked energy services potentially allow fairer revenues for EV owners where the proposed submetering solution extends the service portfolio using smart meters.

Regarding rate plans offered by utilities, the gains obtained from submetering can be enhanced with high on-peak charges or the super off-peak periods proposed by some electricity providers, for example. We have tuned our counterfactual time-varying energy profiles for a flat energy tariff, so that the energy costs in the end would be practically the same. However, in real-life tariff books, the spread between a flat rate and the price levels of a domestic time-of-use proposed could be higher than that captured in our results. We have shown that the concept can bring fair yearly gains varying from $64 to $110 if well managed. However, the exact benefit remains sensitive to the rate plans available from each utility. Regulators should incentivise utilities to propose more time-varying rate plans to boost demand-side management for users, including those with EVs and DERs.

Finally, two different tariffs give users more flexibility to optimise the electricity bill for a particular load like electric vehicles. Moreover, 'type-of-use' tariffs can also be applied for other loads such as heat pumps. In theory, three tariffs driving different appliances in the same household could bring higher electricity bill savings if well optimised. This framework would substantially complexify users' understanding of their bills, which could discourage them from adopting this solution. Another drawback is the risk of badly designing a rate plan that could over-incentivise a specific type of load usage and bring a deficit to the utility budget. The concept of price discrimination applied in these tariffs should be carefully taken into account to avoid any subsequent fairness issues between customers.

## 7. Conclusion

Adopting an electricity rate designed specifically to charge EVs at home enables a reduction in charging costs for myriad dwellings. Assuming that household demand without flexible loads is quasi-inelastic, separating the billing for EV charging gives a fair incentive to adopt domestic EV-only rates. We have developed a game-theoretical model expressed and treated as a mathematical programme with equilibrium constraints (MPEC) to capture the interaction between a national regulatory authority (NRA) and dwellings. This type of modelling is essential to determine the network tariff endogenously depending on the reactions of all agents. With the uptake of domestic EV-only tariffs, grid operators will also have the necessary information about charging events and could better accommodate other EVs while avoiding costly grid reinforcements. NRAs at their end can adopt a tariff structure that gives incentives for DER and separate ones for EV adoption.

A fair energy cost reduction is observed with an EV-only tariff for the adopter while keeping network charges fixed. However, by recovering grid costs via a three-part tariff that contains capacity charges, the increase in network costs offsets the gains brought by energy savings. With a pure volumetric tariff, fairness issues are nuanced, resulting in well-allocated network costs distribution in which consumers experience a higher decrease in their charges. In addition, submetering can bring yearly gains varying from $64 to $110 with this type of tariff.

Our results could be used to support the creation of new local projects to demonstrate the cost-benefits of this solution applied to specific contexts. The majority of pilot projects currently concentrated in the U.S. may not sufficiently represent the diversity of contexts in which EV-only tariffs could work properly. Extending the number of pilots to other countries is essential to identify all possible barriers. With respect to future work, an extension of the problem to regulators in neighbouring countries could frame the economic spill-over between them. Different EV penetration, DER levels and adoption rates of EV-only tariffs could provide a more reliable framework. Moreover, customer surveys including the option of EV-only tariffs to assess their acceptance are crucial to better design this type of tariffs.

### CRediT authorship contribution statement

**Icaro Silvestre Freitas Gomes:** Conceptualisation, Methodology, Formal analysis, Investigation, Writing – Original Draft, Writing – Review & Editing. **Adam F. Abdin:** Conceptualisation, Methodology, Formal analysis, Investigation, Writing – Review & Editing. **Jakob Punchinger:** Methodology, Investigation, Writing – Review & Editing, Supervision. **Yannick Perez:** Methodology, Investigation, Writing – Review & Editing, Supervision, Funding acquisition.

### Acknowledgements

The authors have no conflicts of interest.

### Appendix A. The detailed mathematical model

An overview of the bi-level model and the associated sets, variables and parameters is described together with the framework to solve it.

**Symbols**

- Indices:
  - *c:* Customers
  - *h:* Time period (hours)
- Sets:
  - $C$: Set of all customers
  - $C_{EV}$: Subset of customers set containing consumers with submetered EV ($\subset C$)
  - $C_{EV/DER}$: Subset of customers set containing prosumers with submetered EV ($\subset C$)
  - $H$: Set of hours
- Upper-level parameters:
  - *egc:* Existing grid capacity [ kW]
  - $I^{DSO}$*:* Annualised investment cost for grid capacity [ €/kW/year]
  - *n:* number of agents [ -]
  - *NM:* Net metering coefficient [ -]
  - $PBE^{P}_{c,h}$: Price of buying energy for agent c at hour h [€/kWh]
  - $PSE^{P}_{c,h}$: Price of selling energy for agent c at hour h [€/kWh]
  - *SunkCosts:* Sunk annualised grid costs, scaled per average consumer [€]
  - *W:* Weight of hour h [ h/h]
- Lower-level parameters
  - $D^{\Delta EV-}_{c,h}$ : EV driving needs of agent c in hour h [ kWh/h]
  - $D_{c,h}$: Electricity demand of agent c in hour h [ kWh/h]
  - $\overline{E}^{EV}_{c,h}$ : Maximum EV state of charge of agent c in hour h [kWh]
  - $\underline{E}^{EV}_{c,h}$: Minimum EV state of charge of agent c in hour h [kWh]
  - $G^{PV}_{c,h}$: Solar resource availability [kW/kWp]
  - $I^{S}$,$I^{PV}$: Annualised DER investment costs for agent c [$/kW/Year]

- $L^S, L^{EV}$: Battery and EV converter losses [%]
- $P_{c,h}^{EV\,dis}$: EV discharging power of agent c in hour h [kW]
- $P_{c,h}^{EV\,ch}$: EV charging power of agent c in hour h [kW]
- $P_c^{ch}$: Power capacity ratio for battery charging for agent c [kW/kWh]
- $P_c^{dis}$: Power capacity ratio for battery discharging for agent c [kW/kWh]
- $R; R^{EV}$: Battery and EV battery self-discharge [%]
- $S^{\%max}, S^{\%min}$: Maximum and minimum battery allowed state of charge level [%]
- $U_c^S, U_c^{PV}$: Resource limits of storage and PV for agent c [kWh,kW]

- Upper-level variables
  - *agc*: Additional grid capacity investment in interconnection [ kW]
  - *cnt*: Capacity network charge [ €/kW]
  - $e_h^{GE}$: EV fleet exports in hour h [ kWh/h]
  - $e_h^G$: Total net load [ kWh/h]
  - $e_h^{GI}$: Total imports in hour h [ kWh/h]
  - *fnt*: Fixed network charge [ €/customer]
  - *vnt*: Volumetric network charge [ €/kWh]
- Lower-level Variables
  - $d_{c,h}^{\Delta EV+}, d_{c,h}^{\Delta EV-}$: EV battery charge/discharge in hour h [kWh/h]
  - $d^{\Delta+}_{c,h}, d^{\Delta}_{c,h}{}^{-}$: Stationary battery charge/discharge in hour h [kWh/h]
  - $exp^P_{c,h}$: Energy exported to grid by agent c in hour h [kWh/h]
  - $exp^L_{c,h}$: Energy exported to agent behind the meter [kWh/h]
  - $g_{c,h}^{PV}$: Solar PV electricity generation by agent c in hour h [kWh/h]
  - $ic_c^{PV}$: Solar PV installed capacity for agent c [kW]
  - $ic_c^S$: Battery storage capacity for agent c [kWh]
  - $imp^P_{c,h}$: Energy imported from grid by agent c in hour h [kWh/h]
  - $imp^L_{c,h}$: Energy imported from agent behind the meter [kWh/h]
  - $p_c$: Measured peak power of agent c [kW]

- $s^{EV}_{c,h}$: EV battery state of charge in hour h [kWh]
- $s_{c,h}$: Stationary battery state of charge in hour h [kWh]

**Upper-level objective function**:

$$Min\ Cost^{NRA} = Cost^{DER} + Cost^{P} + Cost^{N}\ , \quad (A.1)$$

**where**:

$$Cost^{DER} = \sum_{c\ \in\ C} \left( I^{S}_{c} * ic^{S}_{c}\ +\ I^{PV}_{c} * ic^{PV}_{c} \right) \quad (A.2)$$

$$Cost^{P}\ = \sum_{c\ \in\ C} \sum_{h\ \in\ H} W_h * \left( imp^{P}_{c,h} * PBE^{P}_{c,h}\ -\ exp^{P}_{c,h} * PSE^{P}_{c,h} \right) \quad (A.3)$$

$$Cost^{N} = SunkCosts\ +\ I^{DSO}\ *\ agc \quad (A.4)$$

**Subject to**:

1. Coincidental peak: To maintain the linearity of the problem, Eq. (5) which calculates the hourly sum of net imports and net exports is split into Eqs. (A.5) and (A.6). As just one of them will be non-zero, the sum of both (Eq. (A.7)) measures the aggregate hourly demand:

$$e^{GI}_{h}\ \geqslant \sum_{c\ \in\ C} \left( imp^{P}_{c,h}\ -\ exp^{P}_{c,h} \right),\ \forall h\ \in\ H \quad (A.5)$$

$$e^{GE}_{h}\ \geqslant \sum_{c\ \in\ C} \left( exp^{P}_{c,h}\ -\ imp^{P}_{c,h} \right), \forall h\ \in\ H \quad (A.6)$$

$$e^{G}_{h}\ =\ e^{GI}_{h}\ +\ e^{GE}_{h},\ \forall h\ \in\ H \quad (A.7)$$

2. Grid capacity

$$egc + agc > e^{G}{}_{h},\ \forall h \in H \quad \text{(A.8)}$$

3. DSO cost recovery

$$Cost^{N} = \sum_{c\,\in\,C}\sum_{h\,\in\,H} W*\left(imp^{P}_{c,h} - NM * exp^{P}_{c,h}\right) * vnt + \sum_{c\,\in\,C} p_c * cnt + n * fnt \quad \text{(A.9)}$$

**Lower-level objective function**:

$$Min\ Cost = Cost^{DER}_{c} + Cost^{P}_{c} + Cost^{N}_{c}, \quad \text{(A.10)}$$

**where**:

$$Cost^{DER}_{c} = I^{S} * ic^{S}_{c} + I^{PV} * ic^{PV}_{c}$$

$$Cost^{P}_{c} = \sum_{h\,\in\,H} W*\left(imp^{P}_{c,h} * PBE^{P}_{c,h} - exp^{P}_{c,h} * PSE^{P}_{c,h}\right)$$

$$Cost^{N}_{c} = \sum_{h\,\in\,H} W*\left(imp^{P}_{c,h} - NM * exp^{P}_{c,h}\right) * vnt + p_c * cnt + fnt$$

(A.11)

(A.12)

(A.13)

**Subject to**:

1. Energy balance

$$-D_{c,h} - d^{\Delta EV+}_{c,h} + d^{\Delta EV-}_{c,h} - d^{\Delta+}_{c,h} + d^{\Delta-}_{c,h} + ic^{PV}_{c} * G^{PV}_{c,h} + imp^{P}_{c,h} - exp^{P}_{c,h} + \alpha_c * \left(imp^{L}_{c,h} - exp^{L}_{c,h}\right) = 0\ \forall c \in C, h \in H : \left(\lambda^{EB}_{c,h}\right) \quad \text{(A.14)}$$

2. Peak power measurement

$$-p_c + imp^{P}_{c,h} + exp^{P}_{c,h} \leqslant 0\ \forall c \in C, h \in H : \left(\mu^{G}_{c,h}\right); \quad \text{(A.15)}$$

3. EV storage constraints

$$s_{c,h}^{EV} - s_{c,h-1}^{EV} * \left(1 - R^{EV}\right) - d_{c,h}^{\Delta EV+} * \left(1 - L^{EV}\right) + d_{c,h}^{\Delta EV-} * \left(1 + L^{EV}\right) + D_{c,h}^{\Delta EV-} = 0,$$
$$\forall c \in C, h \in H \setminus \{1\} : \left(\lambda_{c,h}^{EV1}\right) \quad \text{(A.16)}$$

$$s_{c,1}^{EV} - SOC_0^{EV} - d_c^{\angle}$$

$$s_{c,H}^{EV} - SOC_0^{EV} = 0,\ \forall c \in C, h \in H : \left(\lambda_c^{EV2}\right)$$
$$s_{c,h}^{EV} - \overline{E}_{c,h}^{EV} \leqslant 0,\ \forall c \in C, h \in H : \left(\mu_{c,h}^{EV2}\right)$$
$$\underline{E_{c,h}^{EV}} - s_{c,h}^{EV} \leqslant 0,\ \forall c \in C, h \in H : \left(\mu_{c,h}^{EV3}\right)$$
$$d_{c,h}^{\Delta EV+} - P_{c,h}^{EV ch} \leqslant 0,\ \forall c \in C, h \in H : \left(\mu_{c,h}^{EV4}\right)$$
$$d_{c,h}^{\Delta EV-} - P_{c,h}^{EV dis} \leqslant 0,\ \forall c \in C, h \in H : \left(\mu_{c,h}^{EV5}\right)$$

(A.18)

(A.19)

(A.20)

(A.21)

(A.22)

4. Battery storage constraints

$$s_{c,h} - s_{c,h-1} * \left(1 - R^{S}\right) - d_{c,h}^{\Delta+} * \left(1 - L^{S}\right) + d_{c,h}^{\Delta-} * \left(1 + L^{S}\right) = 0,\ \forall c \in C, h \in H \setminus \{1\} : \left(\lambda_{c,h}^{S1}\right) \quad \text{(A.23)}$$

$$s_{c,1} - s_{c,H} * \left(1 - R^{S}\right) - d_{1}^{\Delta+} * \left(1 - L^{S}\right) + d_{1}^{\Delta-} * \left(1 + L^{S}\right) = 0,\ \forall c \in C : \left(\lambda_{c,1}^{S1}\right) \quad \text{(A.24)}$$

$$ic_c^{s} - U_c^{S} \leqslant 0,\ \forall c \in C : \left(\mu_c^{S1}\right) \quad \text{(A.25)}$$

$$s_{c,h} - ic_c^{s} * S^{\%max} \leqslant 0\ \forall c \in C, h \in H : \left(\mu_{c,h}^{S2}\right) \quad \text{(A.26)}$$

(A.27)

$$d_{c,h}^{\Delta+} - ic_c^{s} * P_c^{ch} \leqslant 0\ \forall c \in C, h \in H : \left(\mu_{c,h}^{S4}\right) \quad \text{(A.28)}$$

$$d_{c,h}^{\Delta-} - ic_c^{s} * P_c^{dis} \leqslant 0\ \forall c \in C, h \in H : \left(\mu_{c,h}^{S5}\right) \quad \text{(A.29)}$$

5. Solar PV constraints

$$ic_c^{PV} - U_c^{PV} \leqslant 0\ \forall c \in C : \left(\mu_c^{PV1}\right) \quad \text{(A.30)}$$

$$\sum_{c\ \in\ C_{EV}} \left(imp^{L}_{c,h}\ -\ exp^{L}_{c,h}\right) = 0\ ,\forall h\ \in\ H\ :\ \left(\lambda_{h}^{L_{EV}}\right)$$

$$\sum_{c\ \in\ C_{EV/DER}} \left(imp^{L}_{c,h}\ -\ exp^{L}_{c,h}\right) = 0\ ,\forall h\ \in\ H\ :\ \left(\lambda_{h}^{L_{EV/DER}}\right)$$

$$imp^{L}_{EV,h}\ -\ d^{\Delta-}_{c,h}\ -\ ic^{PV}_{c} * G^{PV}_{c,h}\ \leqslant 0\ \forall c\ \in (C_{EV} \cup C_{EV/DER})\ ,\ \forall h\ \in H : \left(\mu^{impL}_{c,h}\right)$$

$$imp^{L}_{c,h}\ -\ d^{\Delta EV-}_{EV,h}\ \leqslant\ 0\ \forall c\ \in C_{EV}, \forall c\ \in (C_{EV} \cup C_{EV/DER})\ ,\ \forall h\ \in H\ : \left(\mu^{impL2}_{c,h}\right)$$

6. Submetering constraints

(A.31)

(A.32)

(A.33)

(A.34)

7. Non-negativity constraints

$$-\,imp^{P}_{c,h} \leqslant 0 \;\forall c \in C, h \in H \;:\; \left(\mu^{imp}_{c,h}\right) \tag{A.35}$$

$$-\,exp^{P}_{c,h} \leqslant 0 \;\forall c \in C, h \in H \;:\; \left(\mu^{exp}_{c,h}\right) \tag{A.36}$$

$$-\,d^{\Delta -}_{c,h} \leqslant 0 \;\forall c \in C, h \in H \;:\; \left(\mu^{dminb}_{c,h}\right) \tag{A.37}$$

$$-\,d^{\Delta +}_{c,h} \leqslant 0 \;\forall c \in C, h \in H \;:\; \left(\mu^{dplusb}_{c,h}\right) \tag{A.38}$$

$$-\,d^{\Delta EV+}_{c,h} \leqslant 0 \;\forall c \in C, h \in H \;:\; \left(\mu^{dEV+}_{c,h}\right) \tag{A.39}$$

$$-\,d^{\Delta EV-}_{c,h} \leqslant 0 \;\forall c \in C, h \in H \;:\; \left(\mu^{dEV-}_{c,h}\right) \tag{A.40}$$

$$-\,ic^{PV}_{c} \leqslant 0 \;\forall c \in C, h \in H \;:\; \left(\mu^{icpv}_{c}\right) \tag{A.41}$$

$$-\,ic^{S}_{c} \leqslant 0 \;\forall c \in C, h \in H \;:\; \left(\mu^{ics}_{c}\right) \tag{A.42}$$

$$-\,s^{EV}_{c,h} \leqslant 0 \;\forall c \in C, h \in H \;:\; \left(\mu^{sEV}_{c,h}\right) \tag{A.43}$$

$$-\,s_{c,h} \leqslant 0 \;\forall c \in C, h \in H \;:\; \left(\mu^{s}_{c,h}\right) \tag{A.44}$$

$$-\,imp^{L}_{c,h} \leqslant 0 \;\forall c \in C, h \in H \;:\; \left(\mu^{Limp}_{c,h}\right) \tag{A.45}$$

$$-\,exp^{L}_{c,h} \leqslant 0 \;\forall c \in C, h \in H \;:\; \left(\mu^{Lexp}_{c,h}\right) \tag{A.46}$$

*Appendix A.1. Transforming the bi-level problem into a solvable MPEC*

The lower-level conditions are changed by their KKT optimality conditions. This step allows transformation of the bi-level problem into a MPEC which now has a single objective function. We derive the KKT necessary optimality conditions from the primal feasibility restrictions of the lower problem:

$$W * \left(PBE_{c,h}^{P} + vnt\right) + \lambda_{c,h}^{EB} + \mu_{c,h}^{G} - \mu_{c,h}^{imp} = 0 , \forall c \in C, h \in H$$

$$-W * \left(PSE_{c,h}^{P} + NM * vnt\right) - \lambda_{c,h}^{EB} + \mu_{c,h}^{G} - \mu_{c,h}^{exp} = 0 , \forall c \in C, h \in H$$

$$I^{PV} + \sum_{h \in H} \lambda_{c,h}^{EB} * G_{c,h}^{PV} + \mu_{c}^{PV1} - \mu_{c}^{icpv} = 0 , \forall c \in C$$

$$-\lambda_{c,h}^{EB} - \left(1 - L^{EV}\right) * \lambda_{c,h}^{EV1} + \mu_{c,h}^{EV4} - \mu_{c,h}^{dEV+} = 0 , \forall c \in C, h \in H$$

$$\lambda_{c,h}^{EB} + \left(1 + L^{EV}\right) * \lambda_{c,h}^{EV1} + \mu_{c,h}^{EV5} - \mu_{c,h}^{dEV-} = 0 , \forall c \in C, h \in H$$

$$\lambda_{c,h}^{EV1} - \left(1 - R^{EV}\right) * \lambda_{c,h+1}^{EV1} + \mu_{c,h}^{EV2} - \mu_{c,h}^{EV3} - \mu_{c,h}^{sEV} = 0 , \forall c \in C, h \in H \setminus \{48\}$$

$$\lambda_{c,H}^{EV1} + \lambda_{c}^{EV2} + \mu_{c,H}^{EV2} - \mu_{c,H}^{EV3} - \mu_{c,H}^{sEV} = 0 \ \forall c \in C, h = H$$

$$-\lambda_{c,h}^{EB} - \left(1 - L^{S}\right) * \lambda_{c,h}^{S1} + \mu_{c,h}^{S4} - \mu_{c,h}^{d+} = 0 , \forall c \in C, h \in H$$

$$\lambda_{c,h}^{EB} + \left(1 + L^{S}\right) * \lambda_{c,h}^{S1} + \mu_{c,h}^{S5} - \mu_{c,h}^{d-} = 0 , \forall c \in C, h \in H$$

$$\lambda_{c,h}^{S1} - \left(1 - R^{S}\right) * \lambda_{c,h+1}^{S1} + \mu_{c,h}^{S2} - \mu_{c,h}^{S3} - \mu_{c,h}^{s} = 0 , \forall c \in C, h \in H \setminus \{48\}$$

$$\lambda_{c,H}^{S1} - \left(1 - R^{S}\right) * \lambda_{c,1}^{S1} + \mu_{c,H}^{S2} - \mu_{c,H}^{S3} - \mu_{c,H}^{s} = 0 , \forall c \in C$$

$$cnt - \sum_{h \in H} \mu_{c,h}^{G} = 0 , \forall c \in C$$

$$\lambda_{c,h}^{EB} + \lambda_{h}^{L_{EV}} - \mu_{c,h}^{impL} = 0 , \forall c \in C, h \in H$$

$$-\lambda_{c,h}^{EB} - \lambda_{h}^{L_{EV}} - \mu_{c,h}^{expL} = 0 , \forall c \in C, h \in H$$

(A.47)

(A.48)

(A.49)

(A.50)

(A.51)

(A.52)

(A.53)

(A.54)

(A.55)

(A.56)

(A.57)

(A.58)

(A.59)

(A.60) $$0 \leq imp^{P}_{c,h} \perp \mu^{imp}_{c,h} \geqslant 0, \forall c \in C, h \in H$$

(A.61)

$$0 \leq exp^{P}_{c,h} \perp \mu^{exp}_{c,h} \geqslant 0, \forall c \in C, h \in H \quad (A.62)$$

$$0 \leq ic^{PV}_{c} \perp \mu^{icpv}_{c,h} \geqslant 0, \forall c \in C \quad (A.63)$$

$$0 \leq U^{PV}_{c} - ic^{PV}_{c} \perp \mu^{PV1}_{c} \geqslant 0, \forall c \in C \quad (A.64)$$

(A.65)

$$0 \leq d^{\Delta EV-}_{c,h} \perp \mu^{dEV-}_{c,h} \geqslant 0, \forall c \in C, h \in H \quad (A.66)$$

$$0 \leq s^{EV}_{c,h} \perp \mu^{sEV}_{c,h} \geqslant 0, \forall c \in C, h \in H \quad (A.67)$$

$$0 \leq \overline{E}^{EV}_{c,h} - s^{EV}_{c,h} \perp \mu^{EV2}_{c,h} \geqslant 0, \forall c \in C, h \in H$$

$$0 \leq s^{EV}_{c,h} - \underline{E}^{EV}_{c,h} \perp \mu^{EV3}_{c,h} \geqslant 0, \forall c \in C, h \in H \quad (A.68)$$

$$0 \leq P^{EVch}_{c,h} - d^{\Delta EV+}_{c,h} \perp \mu^{EV4}_{c,h} \geqslant 0, \forall c \in C, h \in H \quad (A.69)$$

(A.70)

$$0 \leq P^{EVdis}_{c,h} - d^{\Delta EV-}_{c,h} \perp \mu^{EV5}_{c,h} \geqslant 0, \forall c \in C, h \in H \quad (A.71)$$

$$0 \leq d^{\Delta+}_{c,h} \perp \mu^{d+}_{c,h} \geqslant 0, \forall c \in C, h \in H \quad (A.72)$$

$$0 \leq d^{\Delta-}_{c,h} \perp \mu^{d-}_{c,h} \geqslant 0, \forall c \in C, h \in H \quad (A.73)$$

$$0 \leq s_{c,h} \perp \mu^{s}_{c,h} \geqslant 0, \forall c \in C, h \in H \quad (A.74)$$

$$0 \leq ic^{s}_{c} * S^{\%max} - s_{c,h} \perp \mu^{S2}_{c,h} \geqslant 0, \forall c \in C, h \in H \quad (A.75)$$

$$0 \leq s_{c,h} - S^{\%min} * ic^{s}_{c} \perp \mu^{S3}_{c,h} \geqslant 0, \forall c \in C, h \in H \quad (A.76)$$

$$0 \leq ic^{s}_{c} * P^{ch}_{c} - d^{\Delta+}_{c,h} \perp \mu^{S4}_{c,h} \geqslant 0, \forall c \in C, h \in H \quad (A.77)$$

$$0 \leq ic^{s}_{c} * P^{dis}_{c} - d^{\Delta-}_{c,h} \perp \mu^{S5}_{c,h} \geqslant 0, \forall c \in C, h \in H \quad (A.78)$$

$$0 \leq ic^{s}_{c} * P^{ch}_{c} - d^{\Delta+}_{c,h} \perp \mu^{S4}_{c,h} \geqslant 0, \forall c \in C, h \in H \quad (A.79)$$

$$0 \leq imp^{L}_{c,h} \perp \mu^{impL}_{c,h} \geqslant 0, \forall c \in C, h \in H \quad (A.80)$$

$$0 \leq exp^{L}_{c,h} \perp \mu^{expL}_{c,h} \geqslant 0, \forall c \in C, h \in H \quad (A.81)$$

$$0 \leqslant -imp^{L}_{EV,h} + d^{\Delta EV-}_{c,h} + d^{\Delta-}_{c,h} + ic^{PV}_{c} * G^{PV}_{c,h} \perp \mu^{impL1}_{c,h} \geqslant 0 \;\; \forall c \in (C_{EV} \cup C_{EV/DER}), \forall h \in$$

(A.82)

$$0 \leqslant -imp^{L}_{c,h} + d^{\Delta EV-}_{EV,h} \perp \mu^{impL2}_{c,h} \geqslant 0 \;\; \forall c \in (C_{EV} \cup C_{EV/DER}), \forall h \in H \quad (A.83)$$

$$-D_{c,h} - d^{\Delta EV+}_{c,h} + d^{\Delta EV-}_{c,h} - d^{\Delta+}_{c,h} + d^{\Delta-}_{c,h} + ic^{PV}_{c} * G^{PV}_{c,h}$$

$$imp^{P}_{c,h} - exp^{P}_{c,h} + \alpha_c * (imp^{L}_{c,h} - exp^{L}_{c,h}) = 0 \; : \; (\lambda^{EB}_{c,h}, free), \forall c \in C, h \in H \quad (A.84)$$

$$s_{c,h}^{EV} - s_{c,h-1}^{EV} * \left(1 - R^{EV}\right) - d_{c,h}^{\Delta EV+} * \left(1 - L^{EV}\right) + d_{c,h}^{\Delta EV-} * \left(1 + L^{EV}\right) + D_{c,h}^{\Delta EV-} = 0 \ : \left(\lambda_{c,h}^{EV1}, free\right) \quad \forall c \in C, h \in H \setminus \{1\} \quad (A.85)$$

$$s_{c,1}^{EV} - SOC_0^{EV} - d_{c,1}^{\Delta EV+} * \left(1 - L^{EV}\right) + d_{c,1}^{\Delta EV-} * \left(1 + L^{EV}\right) + D_{c,1}^{\Delta EV-} = 0 \ : \left(\lambda_{c,1}^{EV1}, free\right) \ , \forall c \in C \quad (A.86)$$

$$s_{c,H}^{EV} - SOC_0^{EV} = 0 \ : \left(\lambda_c^{EV2}, free\right) \ , \forall c \in C \quad (A.87)$$

$$s_{c,h} - s_{c,h-1} * \left(1 - R^S\right) - d_{c,h}^{\Delta+} * \left(1 - L^S\right) + d_{c,h}^{\Delta-} * \left(1 + L^S\right) = 0 \ : \left(\lambda_{c,h}^{S1}, free\right), \forall c \in C, h \in H \setminus \{1\} \quad (A.88)$$

$$s_{c,1} - s_{c,H} * \left(1 - R^S\right) - d_{c,1}^{\Delta+} * \left(1 - L^S\right) + d_{c,1}^{\Delta-} * \left(1 + L^S\right) = 0 \ : \left(\lambda_{c,1}^{S1}, free\right), \forall c \in C \quad (A.89)$$

$$\sum_{c \in C_{EV}} \left(imp_{c,h}^L - exp_{c,h}^L\right) = 0 \ : \left(\lambda_h^{L_{EV}}, free\right) \ , \forall h \in H \quad (A.90)$$

$$\sum_{c \in C_{EV/DER}} \left(imp_{c,h}^L - exp_{c,h}^L\right) = 0 \ : \left(\lambda_h^{L_{EV/DER}}, free\right) \ , h \in H \quad (A.91)$$

Before treating the non-linearities of the model, an adjustment regarding cost recovery equality is needed to facilitate convergence. As in Schittekatte and Meeus (2020), the total network charge costs collected should be within a band (calibrated as $\delta = 0.1\%$):

$$\left(Cost^N\right) * (1-\delta) - \sum_{c \in C} \sum_{h \in H} W * \left(imp_{c,h}^P - NM * exp_{c,h}^P\right) * vnt + \sum_{c \in C} p_c * cnt + n * fnt) \leqslant 0 \quad (A.92)$$

$$-\left(Cost^N\right) * (1+\delta) + \sum_{c \in C} \sum_{h \in H} W * \left(imp_{c,h}^P - NM * exp_{c,h}^P\right) * vnt + \sum_{c \in C} p_c * cnt + n * fnt \leqslant 0 \quad (A.93)$$

To fully transform the MPEC into a MILP, two non-linearities must be taken into account: the bilinear terms and the complementarity constraints in the KKT conditions. First, the bilinear products in the equality constraint of the upper-level ($imp_{c,h}^P * vnt$, $exp_{c,h}^P * vnt$ and $p_c * cnt$) are already taken into account internally by solver Gurobi 9.1 (Gurobi-Optimisation, 2021). Instead of discretising the terms by using binary expansion (as in Momber (2015), page 102) beforehand, the solver deals with this type of non-convexity by applying cutting planes and special branching techniques. Finally, the complementarity constraints are linearised using the Fotuny-Amat method (Fortuny-Amat et al. (1981)) in which they are reformulated using additional binary variables and large enough constants (Big-Ms). Another solution to deal with the complementarity constraints in the MPEC framework could be SOS1 variables as proposed in Siddiqui and Gabriel (2013). At this point, the bi-level problem turns into a mixed-integer linear programme that can be solved using the prior Gurobi 9.1 solver:

$$imp_{c,h}^{P} \le M^{imp} * \left(1 - r_{c,h}^{a}\right) \quad ,\forall c \in C, h \in H$$

$$W * \left(PBE_{c,h}^{P} + vnt\right) + \lambda_{c,h}^{EB} + \mu_{c,h}^{G} \le M^{imp} * r_{c,h}^{a} \; \forall c \in C, h \in H$$

$$exp_{c,h}^{P} \le M^{exp} * \left(1 - r_{c,h}^{b}\right) \quad ,\forall c \in C, h \in H$$

$$-W * \left(PSE_{c,h}^{P} + NM * vnt\right) - \lambda_{c,h}^{EB} + \mu_{c,h}^{G} \le M^{exp} * r_{c,h}^{b} \; ,\forall c \in C, h \in H$$

$$ic_{c}^{PV} \le M^{icpv} * \left(1 - r_{c}^{c}\right) \quad ,\forall c \in C$$

$$I^{PV} + \sum_{h \in H} \lambda_{c,h}^{EB} * G_{c,h}^{PV} + \mu_{c}^{PV1} \le M^{icpv} * r_{c}^{c} \; \forall c \in C$$

$$U_{c}^{PV} - ic_{c}^{PV} \le M^{PV1} * \left(1 - r_{c}^{d}\right) \quad ,\forall c \in C$$

$$\mu_{c}^{PV1} \le M^{PV1} * r_{c}^{d} \; ,\forall c \in C$$

$$d_{c,h}^{\Delta EV+} \le M^{dEV+} * \left(1 - r_{c,h}^{e}\right) \; ,\forall c \in C, h \in H$$

$$-\lambda_{c,h}^{EB} - \left(1 - L^{EV}\right) * \lambda_{c,h}^{EV1} + \mu_{c,h}^{EV4} \le M^{dEV+} * r_{c,h}^{e} \; ,\forall c \in C, h \in H$$

$$d_{c,h}^{\Delta EV-} \le M^{dEV-} * \left(1 - r_{c,h}^{f}\right) \quad ,\forall c \in C, h \in H$$

$$\lambda_{c,h}^{EB} + \left(1 + L^{EV}\right) * \lambda_{c,h}^{EV1} + \mu_{c,h}^{EV5} \le M^{dEV-} * r_{c,h}^{f} \quad ,\forall c \in C, h \in H$$

$$s_{c,h}^{EV} \le M^{sEV} * \left(1 - r_{c,h}^{g}\right) \quad ,\forall c \in C, h \in H$$

$$\mu_{c,h}^{sEV} \le M^{sEV} * r_{c,h}^{g} \; ,\forall c \in C, h \in H$$

$$\overline{E}_{c,h}^{EV} - s_{c,h}^{EV} \le M^{uEV2} * \left(1 - r_{c,h}^{h}\right) \quad ,\forall c \in C, h \in H$$

$$\mu_{c,h}^{EV2} \le M^{uEV2} * r_{c,h}^{h} \; ,\forall c \in C, h \in H$$

$$s_{c,h}^{EV} - \underline{E}_{c,h}^{EV} \le M^{uEV3} * \left(1 - r_{c,h}^{i}\right) \quad ,\forall c \in C, h \in H$$

$$\mu_{c,h}^{EV3} \le M^{uEV3} * r_{c,h}^{i} \; ,\forall c \in C, h \in H$$

$$P_{c,h}^{EVch} - d_{c,h}^{\Delta EV+} \le M^{uEV4} * \left(1 - r_{c,h}^{j}\right) \quad ,\forall c \in C, h \in H$$

$$\mu_{c,h}^{EV4} \le M^{uEV4} * r_{c,h}^{j} \; ,\forall c \in C, h \in H$$

$$P_{c,h}^{EVdis} - d_{c,h}^{\Delta EV-} \le M^{uEV5} * \left(1 - r_{c,h}^{k}\right) \quad ,\forall c \in C, h \in H$$

$$\mu_{c,h}^{EV5} \le M^{uEV5} * r_{c,h}^{k} \; ,\forall c \in C, h \in H$$

(A.94)

(A.95)

(A.96)

(A.97)

(A.98)

(A.99)

(A.100)

(A.101)

(A.102)

(A.103)

(A.104)

(A.105)

(A.106)

(A.107)

(A.108)

(A.109)

(A.110)

(A.111) (A.112)

(A.113)

(A.114) (A.115) (A.116)

(A.117)

(A.118)

(A.119)

(A.120)

(A.121)

(A.122)

(A.123)

(A.124)

(A.125) (A.126)

(A.127)

$$\mu^{S2}_{c,h} \leq M^{uS2} * r^{o}_{c,h} \ ,\forall c \in C, h \in H$$

$$s_{c,h} - S^{\%min} * ic^{s}_{c} \leq M^{uS3} * \left(1 - r^{p}_{c,h}\right) \ \forall c \in C, h \in H$$

$$ic^{s}_{c} * S^{\%max} - s_{c,h} \leq M$$

$$\mu^{S3}_{c,h} \leq M^{uS3} * r^{p}_{c,h} \ ,\forall c \in C, h \in H$$

$$ic^{s}_{c} * P^{ch}_{c} - d^{\Delta+}_{c,h} \leq M^{uS4} * \left(1 - r^{q}_{c,h}\right) \ ,\forall c \in C, h \in H$$

$$\mu^{S4}_{c,h} \leq M^{uS4} * r^{q}_{c,h} \ ,\forall c \in C, h \in H$$

$$ic^{s}_{c} * P^{dis}_{c} - d^{\Delta-}_{c,h} \leq M^{uS5} * \left(1 - r^{r}_{c,h}\right) \ ,\forall c \in C, h \in H$$

$$\mu^{S5}_{c,h} \leq M^{uS5} * r^{r}_{c,h} \ ,\forall c \in C, h \in H$$

$$d^{\Delta+}_{c,h} \leq M^{d+} * \left(1 - r^{l}_{c,h}\right) \ ,\forall c \in C, h \in H$$

$$-\lambda^{EB}_{c,h} - \left(1 - L^{S}\right) * \lambda^{S1}_{c,h} + \mu^{S4}_{c,h} \leq M^{d+} * r^{l}_{c,h} \ ,\forall c \in C, h \in H$$

$$d^{\Delta-}_{c,h} \leq M^{d-} * \left(1 - r^{m}_{c,h}\right) \ ,\forall c \in C, h \in H$$

$$\lambda^{EB}_{c,h} + \left(1 + L^{S}\right) * \lambda^{S1}_{c,h} + \mu^{S5}_{c,h} \leq M^{d-} * r^{m}_{c,h} \ ,\forall c \in C, h \in H$$

$$s_{c,h} \leq M^{s} * \left(1 - r^{n}_{c,h}\right) \ ,\forall c \in C, h \in H$$

$$\mu^{s}_{c,h} \leq M^{s} * r^{n}_{c,h} \ ,\forall c \in C, h \in H$$

$$\mu^{P}_{c,h} \leq M^{P} * r^{t}_{c,h} \ ,\forall c \in C$$

$$p_c - imp_{c,h} - exp_{c,h} \leq M^{P} * \left(1 - r^{t}_{c,h}\right) ,\forall c \in C$$

$$imp_{c,h} \leq M^{impL} * \left(1 - r^{u}_{c,h}\right) ,\forall c \in C, h \in H$$

$$\lambda^{EB}_{c,h} + \lambda^{L}_{h} + \mu^{impL}_{c,h} \leq M^{impL} * r^{u}_{c,h} \ ,\forall c \in C, h \in H$$

$$exp_{c,h} \leq M^{expL} * \left(1 - r^{v}_{c,h}\right) ,\forall c \in C, h \in H$$

$$-\lambda^{EB}_{c,h} - \lambda^{L}_{h} + \mu^{expL}_{c,h} \leq M^{expL} * r^{v}_{c,h} \ ,\forall c \in C, h \in H$$

$$\mu^{impL1}_{c,h} \leq M^{impL1} * r^{w}_{c,h} \ ,\forall c \in C, h \in H$$

$$-imp^{L}_{EV,h} + d^{\Delta EV-}_{c,h} + d^{\Delta-}_{c,h} + ic^{PV}_{c} * G^{PV}_{c,h} \leq M^{impL1} * \left(1 - r^{w}_{c,h}\right) \ ,\forall c \in C, h \in H$$

$$\mu^{impL2}_{c,h} \leq M^{mpL2} * r^{x}_{c,h} \ \forall c \in (C_{EV} \cup C_{EV/DER}), \forall h \in H$$

$$-imp^{L}_{c,h} + d^{\Delta EV-}_{EV,h} \leq M^{impL2} * \left(1 - r^{x}_{c,h}\right) \ \forall c \in (C_{EV} \cup C_{EV/DER}), \forall h \in H$$

(A.128)

(A.129)

(A.130)

(A.131) (A.132) (A.133) (A.134) (A.135) (A.136)

(A.137)

(A.138)

(A.139)

**Appendix B. Electric Nation Customer Trial: Insights from the pilot**

The Electric Nation pilot was conducted in the U.K. by the electricity distribution network operator for the Midlands, South Wales and the South West called Western Power Distribution. The aim of the project was to better understand the impact of charging at home on electricity distribution networks. The trials happened between January 2017 and 2018 with a total of 673 smart chargers installed in the participants' homes.

From the final database containing more than 157,520 rows of charging episodes and various columns with detailed parameters such as start and stop charging times, kWh consumed, the power level of the EVSE, battery capacity etc., we derive some important parameters to use in our study case, as for example the connection hours. We calculate the disconnection hour based on two observations: the hour in which the amount of energy needed before leaving home is the highest and the greatest difference of available power between two consecutive hours. We calculate the average daily profile of energy consumption per hour before disconnection and the average available power of the fleet.

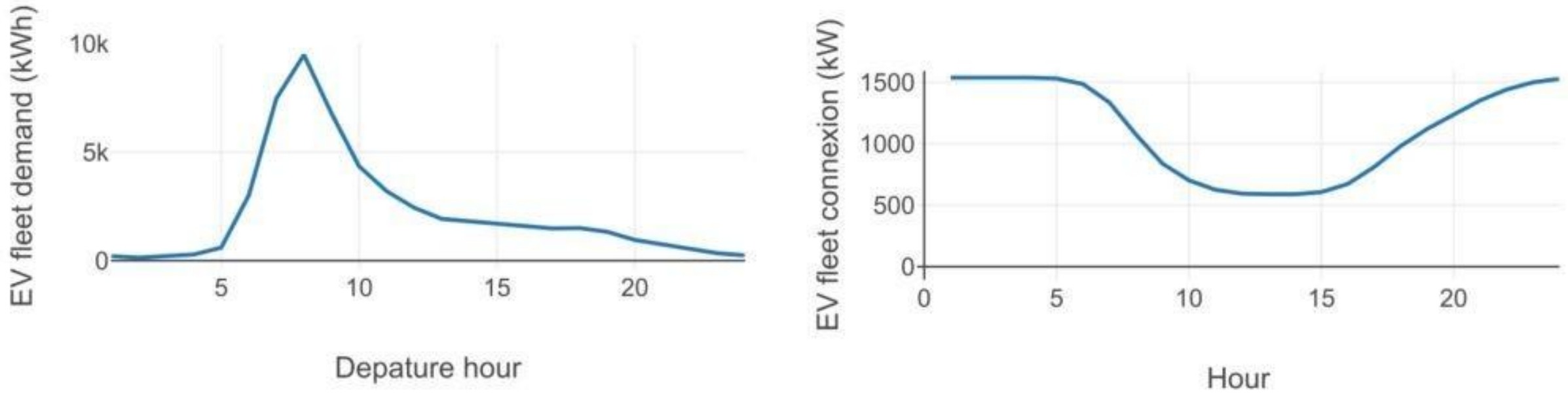

Figure B.10: Electric Nation Trial data.

First, it is observed that the EVs need more energy before leaving at 7am, as is shown in Fig. B.10. Then, the greatest steepness of the curve is between 7am and 8am, meaning that most vehicles left their homes at the former hour. The same analysis is done regarding the reconnection hour, in which the highest power difference between two consecutive hours is from 5pm to 6pm. As a result, the second hour is defined as 5pm.

## Appendix C. Complementary Tables

**Table C.2: Total cost and DER for different energy profiles (MILP analysis)**

| Passive consumers with EV (No V2G) | | | | | |
|---|---|---|---|---|---|
| | One Meter | | | Submetering | |
| | Flat | TOU1 | TOU2 | Flat/TOU1 | Flat/TOU2 |
| Total costs ($) | 1175 | 1050 | 1090 | 1039 | 1079 |
| Coincident peak (kW) | 8,4 | 3,2 | 3,2 | 3,2 | 3,2 |
| **Active consumers with EV (V2G allowed)** | | | | | |
| | One Meter | | | Submetering | |
| | Flat | TOU1 | TOU2 | Flat/TOU1 | Flat/TOU2 |
| Total costs ($) | 1704 | 1323 | 1506 | 1425 | 1533 |
| Coincident peak (kW) | 4,8 | 0 | 0 | 0 | 0 |
| PV investments (kW) | 0,77 | 0,44 | 1,04 | 0,63 | 0,77 |
| Battery investments (kWh) | 0 | 0 | 3,51 | 1,37 | 0 |
| **Active consumers with EV (No V2G)** | | | | | |
| | One Meter | | | Submetering | |
| | Flat | TOU1 | TOU2 | Flat/TOU1 | Flat/TOU2 |
| Total costs ($) | 1704 | 1476 | 1583 | 1568 | 1608 |
| Coincident peak (kW) | 3,1 | 0 | 0 | 3,1 | 3,1 |
| PV investments (kW) | 0,77 | 2,96 | 2,96 | 0,77 | 0,77 |
| Battery investments (kWh) | 0 | 13,6 | 13,6 | 0 | 0 |

Table C.3: Total cost and DER for different energy profiles (PV cost = 1,500$/kWp)

| Passive consumers with EV (No V2G) | | | | | |
|---|---|---|---|---|---|
| | One Meter | | | Submetering | |
| | Flat | TOU1 | TOU2 | Flat/TOU1 | Flat/TOU2 |
| Total costs ($) | 1175 | 1050 | 1090 | 1039 | 1079 |
| Coincident peak (kW) | 8,4 | 3,2 | 3,2 | 3,2 | 3,2 |
| **Active consumers with EV (V2G allowed)** | | | | | |
| | One Meter | | | Submetering | |
| | Flat | TOU1 | TOU2 | Flat/TOU1 | Flat/TOU2 |
| Total costs ($) | 1722 | 1340 | 1529 | 1442 | 1551 |
| Coincident peak (kW) | 9,5 | 0 | 2,5 | 0 | 4,5 |
| PV investments (kW) | 0,45 | 0,36 | 0,45 | 0,38 | 0,45 |
| Battery investments (kWh) | 0 | 0 | 0 | 0 | 0 |
| **Active consumers with EV (No V2G)** | | | | | |
| | One Meter | | | Submetering | |
| | Flat | TOU1 | TOU2 | Flat/TOU1 | Flat/TOU2 |
| Total costs ($) | 1722 | 1598 | 1648 | 1586 | 1627 |
| Coincident peak (kW) | 4,7 | 4,3 | 4,7 | 4,7 | 4,7 |
| PV investments (kW) | 0,45 | 2,5 | 0,77 | 0,45 | 0,45 |
| Battery investments (kWh) | 0 | 11 | 0 | 0 | 0 |

Table C.4: Total network costs for all scenarios according to grid cost structure

| | Sunk | 50/50 | Prospective |
|---|---|---|---|
| Volumetric | 2040 | 1538 | 817 |
| Three-part tariff | 2020 | 1533 | 848 |
| Volumetric (V2G and Submetering) | 2039 | 1618 | 1103 |
| Three-part tariff (V2G and Submetering) | 2035 | 1639 | 1041 |

## References

ACER. Annual Report on the Results of Monitoring the Internal Markets in 2019 Electricity and Natural Gas. 2020.

M. Ansarin, Y. Ghiassi-Farrokhfal, W. Ketter and J. Collins. Cross-subsidies among residential electricity prosumers from tariff design and metering infrastructure. Energy Policy, 145(July):111736, 2020. ISSN 03014215. doi: 10.1016/j.enpol.2020.111736. URL https://doi.org/10.1016/j.enpol.2020.111736.

M. Askeland, S. Backe, S. Bjarghov and M. Korp˚as. Helping end-users help each other: Coordinating development and operation of distributed resources through local power markets and grid tariffs. Energy Economics, 94:105065, 2021. ISSN 0140-9883. doi: 10.1016/j.eneco.2020.105065. URL https://doi.org/10.1016/j.eneco.2020.105065.

M. Avau, N. Govaerts and E. Delarue. Impact of distribution tariffs on prosumer demand response. Energy Policy, 151(February):112116, 2021. ISSN 03014215. doi: 10.1016/j.enpol.2020.112116. URL https://doi.org/10.1016/j.enpol.2020.112116.

S. Backe, G. Kara and A. Tomasgard. Comparing individual and coordinated demand response with dynamic and static power grid tariffs. Energy, 201:117619, 2020. ISSN 0360-5442. doi: 10.1016/j.energy.2020.117619. URL https://doi.org/10.1016/j.energy.2020.117619.

S. Borenstein. The redistributional impact of nonlinear electricity pricing. American Economic Journal: Economic Policy, 4 (3):56–90, 2012. ISSN 19457731. doi: 10.1257/pol.4.3.56.

S. Borenstein, M. Fowlie and J. Sallee. Designing electricity rates for an equitable energy transition. 2021.

S. P. Burger, C. R. Knittel, I. J. P´erez-Arriaga and I. Schneider. The Efficiency and Distributional Effects of Alternative Residential Electricity Rate Designs. The Energy Journal, 41(1):199–240, 2020.

California Energy Commission. Joint iou electric vehicle load research. 2019.

A. Cooper and M. Shuster. Electric Company Smart Meter Deployments: Foundation for a Smart Grid (2021 Update). (April), 2021.

CPUC. Summary of cpuc actions to support zero-emission vehicle adoption. 2021. URL https://www.cpuc.ca.gov/zev/. Accessed: March 25, 2021.

EIA. Electric power monthly. energy information administration. 2019. URL *https: //www.eia.gov/electricity/*. Accessed: March 25, 2021.

EIA. Electricity use in homes. energy information administration. 2020. URL *https: //www.eia.gov/tools/faqs/faq.php*?*id* = 97&*t* = 3. Accessed: June 25, 2021.

Electric Vehicle Database. Energy consumption of full electric vehicles. 2021a. URL *https* : *//ev* − *database.org/cheatsheet/energy* − *consumption* − *electric* − *carl*. Accessed: March 25, 2021.

Electric Vehicle Database. Useable battery capacity of full electric vehicles. 2021b. URL *https: //ev* − *database.org/cheatsheet/useable* − *battery* − *capacity* − *electric* − *car*. Accessed: March 25, 2021.

European Commission. Benchmarking Smart Metering Deployment in EU-28. Number December. 2020. ISBN 9789276172956.

A. J. Fortuny-Amat, B. Mccarl, J. Fortuny-Amat and B. Mccarl. A Representation and Economic Interpretation of a Two-Level Programming Problem. The journal of the Operational Research Society, 32(9):783–792, 1981.

I. Freitas Gomes, Y. Perez and E. Suomalainen. Rate design with distributed energy resources and electric vehicles: A Californian case study. SSRN, (March), 2021. URL https://fsr.eui.eu/publications/?handle=1814/69861.

N. Govaerts, K. Bruninx, H. Le Cadre, L. Meeus and E. Delarue. Spillover effects of distribution grid tariffs in the internal electricity market: An argument for harmonization? Energy Economics, 84, 2019. ISSN 01409883. doi: 10.1016/j.eneco.2019.07.019.

Gurobi-Optimization. Version 9.1: Documentation. 2021. URL *https: //www.gurobi.com/documentation/*9.1*/examples/ bilinear$_p$y.html*.

R. Hemmati and H. Mehrjerdi. Investment deferral by optimal utilizing vehicle to grid in solar powered active distribution networks. Journal of Energy Storage, 30(April):101512, 2020. ISSN 2352152X. doi: 10.1016/j.est.2020.101512. URL https://doi.org/10.1016/j.est.2020.101512.

J. Hildermeier, C. Kolokathis, J. Rosenow, M. Hogan, C. Wiese and A. Jahn. Smart EV charging: A global review of promising practices. World Electric Vehicle Journal, 10(4):1–13, 2019. ISSN 20326653. doi: 10.3390/wevj10040080.

Q. Hoarau and Y. Perez. Interactions between electric mobility and photovoltaic generation: A review. Renewable and Sustainable Energy Reviews, 94(October):510–522, 2018. ISSN 18790690. doi: 10.1016/j.rser.2018.06.039. URL https://doi.org/10.1016/j.rser.2018.06.039.

Q. Hoarau and Y. Perez. Network tariff design with prosumers and electromobility: Who wins, who loses? Energy Economics, 83:26–39, 2019. ISSN 0140-9883. doi: 10.1016/j.eneco.2019.05.009. URL https://doi.org/10.1016/j.eneco.2019.05.009.

IEA. Tracking transport report. 2020.

IEA. Global EV Outlook 2021. Global EV Outlook 2021, 2021. doi: 10.1787/d394399e-en.

C. IOUs. Distributed generation successor tariff workshop: Joint iou proposal. 2021. URL https://www.cpuc.ca.gov/WorkArea/DownloadAsset.aspx?id=6442468353, Accessed: May 24, 2021.

C. King and B. Datta. EV charging tariffs that work for EV owners, utilities and society. Electricity Journal, 31(9):24–27, 2018. ISSN 10406190. doi: 10.1016/j.tej.2018.10.010. URL https://doi.org/10.1016/j.tej.2018.10.010.

K. Knezovi´c, M. Marinelli, A. Zecchino, P. B. Andersen and C. Traeholt. Supporting involvement of electric vehicles in distribution grids: Lowering the barriers for a proactive integration. Energy, 134:458–468, 2017. ISSN 03605442. doi: 10.1016/j.energy.2017.06.075.

S. Ku¨feog˘lu, D. A. Melchiorre and K. Kotilainen. Understanding tariff designs and consumer behaviour to employ electric vehicles for secondary purposes in the United Kingdom. Electricity Journal, 32(6):1–6, 2019. ISSN 10406190. doi: 10.1016/j.tej.2019.05.011.

A. Mayol and C. Staropoli. Giving consumers too many choices: a false good idea? A lab experiment on water and electricity tariffs. European Journal of Law and Economics, 51(2):383–410, 2021. ISSN 15729990. doi: 10.1007/s10657-021-09694-6. URL https://doi.org/10.1007/s10657-021-09694-6.

I. Momber. Benefits of Coordinating Plug-In Electric Vehicles in Electric Power Systems. 2015. ISBN 9788460663140. URL http://www.diva-portal.org/smash/record.jsf?pid=diva2:854273.

NREL. Annual Technology Baseline: Commercial PV. (November), 2020. URL https://atb.nrel.gov/electricity/2020/index.php?t=sd.

Odyss´ee Mure. Electricity consumption per dwelling. Odyss´ee-Mure. 2019. URL *https: //www.odyssee − mure.eu/publications/efficiency −by −sector/households/electricity −consumption−dwelling.html*. Accessed: June 25, 2021.

G. Pasaoglu, D. Fiorello, A. Martino, L. Zani, A. Zubaryeva and C. Thiel. Travel patterns and the potential use of electric cars - Results from a direct survey in six European countries. Technological Forecasting and Social Change, 87:51–59, 2014. ISSN 00401625. doi: 10.1016/j.techfore.2013.10.018. URL http://dx.doi.org/10.1016/j.techfore.2013.10.018.

PGE. EV rates. Pacific Gas and Electric Company. 2021. URL *https: //www.pge.com/en$_U$S/residential/rate − plans/rate − plan − options/electric − vehicle − base − plan/electric − vehicle − base − plan.page*. Accessed: March 25, 2021.

F. Salah, J. P. Ilg, C. M. Flath, H. Basse and C. V. Dinther. Impact of electric vehicles on distribution substations: A Swiss case study. Applied Energy, 137:88–96, 2015. ISSN 0306-2619. doi: 10.1016/j.apenergy.2014.09.091. URL http://dx.doi.org/10.1016/j.apenergy.2014.09.091.

SCE. Tariff books. Southern California Edison. 2019. URL https://www.sce.com/regulatory/tariff-books/ rates-pricing-choices. Accessed: July 17, 2019.

T. Schittekatte and L. Meeus. Least-Cost Distribution Network Tariff Design in Theory and Practice. The Energy Journal, 41 (5), 2020. doi: 10.5547/01956574.41.5.tsch.

T. Schittekatte, I. Momber and L. Meeus. Future-proof tariff design: Recovering sunk grid costs in a world where consumers are pushing back. Energy Economics, 70:484–498, 2018. ISSN 01409883. doi: 10.1016/j.eneco.2018.01.028. URL https://doi.org/10.1016/j.eneco.2018.01.028.

SDGE. EV plans. San Diego Gas and Electric. 2021. URL *https: //www.sdge.com/residential/pricing − plans/about − our − pricing − plans/electric − vehicle − plans.* Accessed: March 25, 2021.

S. Siddiqui and S. A. Gabriel. An SOS1-Based Approach for Solving MPECs with a Natural Gas Market Application. Networks and Spatial Economics, 13(2):205–227, 2013. ISSN 1566113X. doi: 10.1007/s11067-012-9178-y.

P. Simshauser and D. Downer. On the Inequity of Flat-rate Electricity Tariffs Author (s): Paul Simshauser and David Downer Source: The Energy Journal, JULY 2016 , Vol . 37 , No . 3 ( JULY 2016 ), pp . 199-229 Published by: International Association for Energy Economics Stable URL: https://www.jstor.org/stable/44075655 On the Inequity of Flat-rate Electricity Tariffs. 37(3):199–229, 2016.

Smart Electric Power Alliance. Residential Electric Vehicle Rates That Work. (November):1–46, 2019.

A. W. Thompson and Y. Perez. Vehicle-to-Everything (V2X ) energy services, value streams and regulatory policy implications. Energy Policy, 11(July):111136, 2019. ISSN 0301-4215. doi: 10.1016/j.enpol.2019.111136. URL https://doi.org/10.1016/j.enpol.2019.111136.

Western Power Distribution. Final Electric Nation Customer Trial Report. 2019.